\documentclass[12pt,a4paper,reqno]{amsart}

\usepackage{amssymb}
\usepackage{amsmath}
\usepackage{latexsym}
\usepackage{exscale}
\usepackage[latin1]{inputenc}
\usepackage{epsfig}
\usepackage{verbatim}
\usepackage{graphicx}
\usepackage{subfigure}

\newcommand{\R}{\mathbb R}

\newcommand{\Z}{\mathbb{Z}}
\newcommand{\N}{\mathbb{N}}

\newcommand{\cA}{\mathcal{A}}  

\newcommand{\la}{\lambda}
\newcommand{\La}{\Lambda}

\newcommand{\Ga}{\Gamma}

\newcommand{\si}{\sigma}
\newcommand{\Si}{\Sigma}

\newcommand{\norm}[1]{\left\Vert#1\right\Vert}
\newcommand{\abs}[1]{\left| #1 \right|}

\numberwithin{equation}{subsection} 

\newtheorem{thm}{Theorem}[subsection] 
\newtheorem{cor}[thm]{Corollary}
\newtheorem{prop}[thm]{Proposition}
\newtheorem{defi}[thm]{Definition}
\newtheorem{rem}[thm]{Remark}

\newtheorem{lem}[thm]{Lemma}

\begin{document}

\allowdisplaybreaks

\title[RESTRICTED NON-LINEAR APPROXIMATION IN
SEQUENCE SPACES]{\large  RESTRICTED NON-LINEAR APPROXIMATION IN
SEQUENCE SPACES AND APPLICATIONS TO WAVELET BASES AND INTERPOLATION}


\author{Eugenio Hernández}

\address{Eugenio Hernández
\\
Departamento de Matem\'aticas
\\
Universidad Aut\'onoma de Ma\-drid
\\
28049 Madrid, Spain}

\email{eugenio.hernandez@uam.es}

\author{Daniel Vera}

\address{Daniel Vera
\\
Departamento de Matem\'aticas
\\
Universidad Aut\'onoma de Madrid
\\
28049 Madrid, Spain}

\email{daniel.vera@uam.es}

\thanks{Research supported by GrantS MTM2007-60952 and MTM2010-16518 of Spain.}

\date{\today}


\maketitle

\begin{abstract}
Restricted non-linear approximation is a type of N-term
approximation where a measure $\nu$ on the index set (rather than
the counting measure) is used to control the number of terms in the
approximation. We show that embeddings for restricted non-linear
approximation spaces in terms of weighted Lorentz sequence spaces
are equivalent to Jackson and Bernstein type inequalities, and also
to the upper and lower Temlyakov property. As applications we obtain
results for wavelet bases in Triebel-Lizorkin spaces by showing the
Temlyakow property in this setting. Moreover, new interpolation
results for Triebel-Lizorkin and Besov spaces are
obtained.\end{abstract}

\section{Introduction}\label{S:intro}
Thresholding of wavelet coefficients is a technique used in image
processing to compress signals or reduce noise. The simplest
thresholding algorithm $T_\varepsilon (\varepsilon>0)$ of a signal
$f$ is obtained by eliminating from a representation of $f$ the
terms whose coefficients have absolute value smaller than
$\varepsilon$.

Although the thresholding approximants $T_\varepsilon(f)$ are
sometimes a visually faithful representation of $f$, they are not
exact, and from a theoretical point of view an error is introduced
if $f$ is replaced by $T_\varepsilon(f)$. Such errors have initially
been measured in the $L^2-$norm, but it is argued in \cite{KP04}
that procedures having small error in $L^p$, or as stated in the
statistical community, small $L^p-$risk, may reflect better the
visual properties of a signal. Observe that in the usual
thresholding the error is measured in the same space as the signal
is represented, usually $L^2$.

A more general situation is considered in \cite{CDH} where the
wavelet coefficients are thresholded from a representation of the
signal in the Hardy space $H^r$, $0<r<\infty$ (recall that $H^r=L^r$
if $1<r<\infty$), but the error is measured in the Hardy space
$H^p$, $0<p<\infty$. They show that this situation is equivalent to
a type of nonlinear approximation, called \textbf{restricted}, in
which a measure $\nu$ on the index set of dyadic cubes of
$\mathbb{R}^d$ is used to control the number of terms in the
approximation. In the classical $n-$term approximation $\nu(Q)=1$,
$Q\in\mathcal{D}$ (counting measure), and in \cite{CDH}
$\nu(Q)=\abs{Q}^{1-p/r}$.

The article \cite{CDH} provides a description of the approximation
spaces in this setting in terms of certain type of discrete
Lorentz spaces, as well as interpolation results for certain pairs
of $H^p$ and Besov spaces. One of the novelties of this article is
that, although the error is measured in $H^p$, the approximation
spaces are not necessarily contained in $H^p$.

The theory of restricted nonlinear approximation was further
developed in \cite{KP06} considering the case of a quasi-Banach
space $\mathbb{X}$, an unconditional basis
$\mathcal{B}=\{\mathbf{e}_I\}_{I\in\mathcal{D}}$, and a measure
$\nu$ on the countable set $\mathcal{D}$. They show that, in this
abstract setting, restricted thresholding and restricted nonlinear
approximation are linked to the $p-$Temlyakov property for $\nu$
(see definition in \cite{KP06}). They also show that this property
is equivalent to certain Jackson and Bernstein type inequalities and
to have the restricted approximation spaces identified as discrete
Lorentz spaces. The approach in \cite{KP06} is that the
approximation spaces are contained in $\mathbb{X}$ and, hence, not
all results in \cite{CDH} can be recovered.

Denote by $S$ the space of all sequences
$\mathbf{s}=\{s_I\}_{I\in\mathcal{D}}$ of complex numbers indexed by
a countable set $\mathcal{D}$. In the present paper we study
restricted nonlinear approximation for quasi-Banach lattices
$f\subset S$ (see definition in section \ref{sS:Seq_Spcs}). Given a
positive measure $\nu$ on $\mathcal{D}$ we define the restricted
approximation spaces $\mathcal{A}^\xi_\mu(f,\nu)$, $0<\xi<\infty$,
$0<\mu\leq\infty$, as subsets of $S$ using $\nu$ to control the
number of terms in the approximation and $f$ to measure the error
(see section \ref{sS:Rstrct_Nonlin_Apprx_in_Seq_Spcs}).

Denote by $\mathcal{E}=\{\mathbf{e}_I\}_{I\in\mathcal{D}}$ the
canonical basis for $S$. We use a weight sequence
$\mathbf{u}=\{u_I\}_{I\in\mathcal{D}}$, $u_I>0$, to control the
weight of each $\mathbf{e}_I$. Discrete Lorentz spaces
$\ell^\mu_\eta(\nu)$ are defined as sequences
$\mathbf{s}=\{s_I\}_{I\in\mathcal{D}}\in S$ using the $\nu$
distribution function of the sequence $\{u_Is_I\}_{I\in\mathcal{D}}$
(see section \ref{sS:Gnral_Dscrt_LorentzSpcs}). Here, $\eta$ is a
function in $\mathbb{W}$ (see section
\ref{sS:weightfuns_discrt_LrntzSpcs}) more general than
$\eta(t)=t^{1/p}$, $0<p<\infty$.

It is shown in subsections \ref{sS:J_type_ineqs} and
\ref{sS:B_type_ineqs} that the condition
\begin{equation}\label{e:intro_charact_of_demofuncts}
C_1\eta_1(\nu(\Gamma)) \leq \norm{\sum_{I\in\Gamma}
\frac{\mathbf{e}_I}{u_I}}_f \leq C_2\eta_2 (\nu(\Gamma))
\end{equation}
for all $\Gamma\subset\mathcal{D}$, with $\nu(\Gamma)<\infty$,
$\eta_1\in \mathbb{W}$ and $\eta_2\in \mathbb{W}_+$, is equivalent
to inclusions between $\mathcal{A}^\xi_\mu(f,\nu)$ and
$\ell^\mu_{t^\xi\eta(t)}({\bf{u}},\nu)$, and also to some Jackson
and Bernstein type inequalities. When $\eta_1(t)=\eta_2(t)=t^{1/p}$,
condition (\ref{e:intro_charact_of_demofuncts}) is called in
\cite{KP06} the $p-$Temlyakov property.

Working with sequence spaces is not a restriction. Lebesgue,
Sobolev, Hardy and Lipschitz spaces all have a sequence space
counterpart when using the $\varphi-$transform (\cite{FJ88},
\cite{FJ90}) or wavelets (\cite{LM86}, \cite{Me90}, \cite{Da92},
\cite{HW96}, \cite{Ma99}, \cite{Al}, \cite{Ky}). More generally, the
Triebel-Lizorkin, $f^s_{p,r}$, and Besov, $b^s_{p,r}$, spaces of
sequences (see section \ref{sS:SeqSpcs_Assoctd_w_SmoothSpcs}) allow
faithful representations of Triebel-Lizorkin,
$F^s_{p,r}(\mathbb{R}^d)$, and Besov, $B^s_{p,r}(\mathbb{R}^d)$,
spaces (these include all the above spaces). When our results are
coupled with the abstract transference framework designed in
\cite{GH04} we recover results for distribution or function spaces,
as the case may be. One reason to consider such general setting,
besides the obvious generalizations, is that measuring the error
$\norm{f-T_\varepsilon (f)}$ in Sobolev spaces, where the smoothing
properties of $f-T_\varepsilon (f)$ are taken into account, may give
a visually more faithful representation of $f$, than when measured
in $L^p$. Observe that two functions may visually be very different
although they may be close in the $L^p$ norm.

In subsection \ref{sS:Rstrct_Apprx_T-L_seq-spcs} we show that
(\ref{e:intro_charact_of_demofuncts}) holds when
$f=f^{s_1}_{p_1,q_1}$ and
$u_I=\norm{\mathbf{e}_I}_{f^{s_2}_{p_2,q_2}}$, with
$\eta_1(t)=\eta_2(t)=t^{1/p_1}$ and $\nu(I)=\abs{I}^\alpha$ if and
only if $\alpha=p_1(\frac{s_2-s_1}{d}-\frac{1}{p_2})+1\not=1$ or if
$\alpha=1$ then $ p_1=q_1$. When the results of subsection
\ref{sS:J_type_ineqs} and \ref{sS:B_type_ineqs} are applied to this
case, we show that restricted approximation spaces of Triebel
-Lizorkin spaces are identified with discrete Lorentz spaces, which
coincide with Besov spaces for some particular values of the
parameters (Lemma \ref{l:LrntzSpcs_r_BsvSpcs_seq}). The results in
\cite{CDH} and \cite{H07} are simple corollaries. We also give a
result about interpolation of Triebel-Lizorkin and Besov spaces
(section \ref{sS:Appl_Real_Interpol}) with less restrictions on the
parameters than those considered in \cite{CDH}.


The organization of the this paper is as follows. Notation,
definitions, results and comments are given in section
\ref{S:Notation_Def_Reslts} which is divided in subsections $2.x$
with $1\leq x\leq 11$. In section \ref{S:Proofs} we prove the
results stated in section \ref{S:Notation_Def_Reslts}. If a
statement of a result is given in subsection $2.x$, its proof can
be found in subsection $3.x$. Be aware that if a subsection $2.x$
only contains notation, definitions and/or comments, but no
statements of results, the corresponding subsection $3.x$ does not
appear in section $3$.

\vskip 1cm
\section{Notations, Definitions, Statements of Results and
Comments}\label{S:Notation_Def_Reslts}
\subsection{Sequence Spaces}\label{sS:Seq_Spcs}
Let $\mathcal{D}$ be a countable (index) set whose elements will
be denoted by $I$. The set $\mathcal{D}$ could be $\N, \Z, ...$
or, as in the applications we have in mind, the countable set of
dyadic cubes on $\R^d$.

Denote by $S=\mathbb{C}^{\mathcal{D}}$ the set of all sequences of
complex numbers $\mathbf{s}=\{s_I\}_{I\in\mathcal{D}}$ defined over
the countable set $\mathcal{D}$. For each $I\in\mathcal{D}$, we
denote by $\mathbf{e}_I$ the element of $S$ with entry $1$ at $I$
and $0$ otherwise. We write
$\mathcal{E}=\{\mathbf{e}_I\}_{I\in\mathcal{D}}$ for the canonical
basis of $S$. We shall use the notation
$\sum_{I\in\Gamma}s_I\mathbf{e}_I, \Gamma\subset\mathcal{D}$, to
denote the element of $S$ whose entry is $s_I$ when $I\in\Gamma$ and
$0$ otherwise. Notice that no meaning of convergence is attached to
the above notation even when $\Gamma$ is not finite.

\begin{defi}\label{d:quasi-Bnch-seq-lattice}
A linear space of sequences $f\subset S$ is a \textbf{quasi-Banach
(sequence) lattice} if there is a quasi-norm $\norm{\cdot}_f$ in
$f$ with respect to which $f$ is complete and satisfies:

(a) Monotonicity: if $\mathbf{t}\in f$ and
$\abs{s_I}\leq\abs{t_I}$ for all $I\in\mathcal{D}$, then
$\mathbf{s}\in f$ and $\norm{\{s_I\}}_f\leq \norm{\{t_I\}}_f$.

(b) If $\mathbf{s}\in f$, then $lim_{n\to \infty} \norm{ s_{I_n}
\mathbf {e}_{I_n}}_f =0$, for some enumeration $\mathfrak{I}=\{I_1,
I_2,\dots\}$.
\end{defi}
We will say that a quasi-Banach (sequence) lattice $f$ is
\textbf{embedded} in $S$, and write $f\hookrightarrow S$ if
\begin{equation}\label{e:def_f-embdd-S}
\lim_{n\rightarrow\infty} \norm{\mathbf{s}^n-\mathbf{s}}_f=0
    \Rightarrow \lim_{n\rightarrow\infty} s_I^{(n)}=s_I \;\;\forall
    I\in\mathcal{D}.
\end{equation}
\begin{rem}\label{r:monotncty_implies_unconditional}
When $\mathcal{E}=\{\mathbf{e}_I\}_{I\in\mathcal{D}}$ is a Schauder
basis for $f$, condition (a) in Definition
\ref{d:quasi-Bnch-seq-lattice} implies that $\mathcal{E}$ is an
unconditional basis for $f$ with constant $C=1$.
\end{rem}

\vskip0.7cm
\subsection{Restricted Non-linear Approximation in Sequence Spaces}
\label{sS:Rstrct_Nonlin_Apprx_in_Seq_Spcs}

In this paper $\nu$ will denote a positive measure on the discrete
set $\mathcal{D}$ such that $\nu(I)> 0$ for all $I\in\mathcal{D}$.
In the classical $N$-term approximation $\nu$ is the counting
measure (i.e. $\nu(I)=1$ for all $I\in\mathcal{D}$), but more
general measures are used in the restricted non-linear approximation
case. The measure $\nu$ will be used to control the number of terms
in the approximation.

\begin{defi}\label{d:Standard_scheme_Rstrct_NonLin_Apprx}
We say that $(f,\nu)$ is a \textbf{standard scheme} (for restricted
non linear approximation) if

i) $f$ is a quasi-Banach (sequence) lattice embedded in $S$.

ii) $\nu$ is a measure on $\mathcal{D}$ as explained in the
   first paragraph in this section.
\end{defi}

Let $(f,\nu)$ be a standard scheme. For $t>0$, define
$$\Sigma_{t,\nu}=\{\mathbf{t}=\sum_{I\in\Gamma}t_I{\mathbf e}_I: \nu(\Gamma)\leq t\}.$$
Notice that $\Sigma_{t,\nu}$ is not linear, but
$\Sigma_{t,\nu}+\Sigma_{t,\nu}\subset\Sigma_{2t,\nu}$.

Given $\mathbf{s}\in S$, the $f$-error (or $f$-risk) of
approximation to $\mathbf{s}$ by elements of $\Sigma_{t,\nu}$ is
given by
$$\sigma_\nu(t,\mathbf{s})=\sigma_\nu(t,\mathbf{s})_f:=\inf_{\mathbf{t}\in\Sigma_{t,\nu}}
    \norm{\mathbf{s}-\mathbf{t}}_f.$$
Notice that elements $\mathbf{s}\in S$ not in $f$ could have
finite $f$-risk since elements of $\Sigma_{t,\nu}$ could have
infinite number of entries.

\begin{defi}\label{d:Rstrct_Apprx_Spcs}
(\textbf{Restricted Approximation Spaces}) Let $(f,\nu)$ be a
standard scheme.

i) For $0<\xi<\infty$ and $0<\mu<\infty$,
$\mathcal{A}^\xi_\mu(f,\nu)$ is defined as the set of all
$\mathbf{s}\in S$ such that
\begin{equation}\label{e:def_RstrctApprxSpcs}
\norm{\mathbf{s}}_{\mathcal{A}^\xi_\mu(f,\nu)}
    :=
    \left(\int_0^\infty[t^\xi\sigma_\nu(t,\mathbf{s})]^\mu\frac{dt}{t}\right)^{1/\mu}<\infty.
\end{equation}

ii) For $0<\xi<\infty$ and $\mu=\infty$,
$\mathcal{A}^\xi_\infty(f,\nu)$ is defined as the set of all
$\mathbf{s}\in S$ such that
\begin{equation}\label{e:def_RstrctApprxSpcs_infty}
\norm{\mathbf{s}}_{\mathcal{A}^\xi_\infty(f,\nu)}
    := \sup_{t>0} t^\xi\sigma_\nu(t,\mathbf{s}) <\infty.
\end{equation}
\end{defi}

Notice that the spaces $\mathcal{A}^\xi_\mu(f,\nu)$ depend on the
canonical basis $\mathcal E$ of $S$. When $f$ are understood, we
will write $\mathcal{A}^\xi_\mu(\nu)$ instead of
$\mathcal{A}^\xi_\mu(f,\nu)$.

\begin{rem}\label{r:replace_def_RstrctApprxSpcs}
If $\mathbf{s}\in f$, using $\sigma_\nu(t,\mathbf{s})\leq
\norm{\mathbf{s}}_f$, it is easy to see that
(\ref{e:def_RstrctApprxSpcs}) can be replaced by
$\norm{\mathbf{s}}_f$ plus the same integral from $1$ to $\infty$.
We need to consider the whole range $0<t<\infty$ since we do not
assume $\mathbf{s}\in f$. Similar remark holds for $\mu=\infty$ in
(\ref{e:def_RstrctApprxSpcs_infty}). Nevertheless, the properties
of the restricted non-linear approximation spaces are the same as
the $N$-term approximation spaces (see \cite{P81} or
\cite{DeL93}).
\end{rem}

By splitting the integral in dyadic pieces and using the
monotonicity of the $f$-error $\sigma_\nu$ we have an equivalent
quasi-norm for the restricted approximation spaces:
\begin{equation}\label{e:equiv_def_RstrctApprxSpcs}
\norm{\mathbf{s}}_{\mathcal{A}^\xi_\mu(\nu)}
    \approx
    \left(\sum_{k=-\infty}^\infty [2^{k\xi}\sigma_\nu(2^k,\mathbf{s})]^\mu\right)^{1/\mu}.
\end{equation}

\vskip0.7cm
\subsection{The Jackson and Bernstein type inequalities}\label{sS:J&B_ineq}
It is well known the fundamental
role played by the Jackson and Bernstein type inequalities in
non-linear approximation theory. Considering our standard scheme
$(f,\nu)$ we give the following definitions.

\begin{defi}\label{d:Jineq_q-B-seq-lattice}
Given $r>0$, a quasi-Banach (sequence) lattice $g\subset S$
satisfies the Jackson's inequality of order $r$ if there exists
$C>0$ such that
$$\sigma_\nu(t,\mathbf{s})\leq Ct^{-r}\norm{\mathbf{s}}_g \;\; \text{for all } \mathbf{s}\in g.$$
\end{defi}
\begin{defi}\label{d:Bineq_q-B-seq-lattice}
Given $r>0$, a quasi-Banach (sequence) lattice $g\subset S$
satisfies the Bernstein's inequality of order $r$ if there exists
$C>0$ such that
$$\norm{\mathbf{t}}_g\leq Ct^r\norm{\mathbf{t}}_f \;\; \text{for all } \mathbf{t}\in\Sigma_{t,\nu}\cap f.$$
\end{defi}
We do not assume in the above definitions that $g\hookrightarrow
f$, but we need to assume $\mathbf{t}\in\Sigma_{t,\nu}\cap f$ for
Definition \ref{d:Bineq_q-B-seq-lattice} to make sense.

\vskip0.7cm
\subsection{Weight functions for discrete Lorentz spaces}\label{sS:weightfuns_discrt_LrntzSpcs}
\begin{defi}\label{d:weightfun_dscrt_LorentzSpcs}
We will denote by $\mathbb{W}$ the set of all continuous functions
$\eta:[0,\infty)\mapsto[0,\infty)$ such that

i) $\eta(0)=0$ and $\lim_{t\rightarrow\infty}\eta(t)=\infty$

ii) $\eta$ is non-decreasing

iii) $\eta$ has the doubling property, that is, there exists $C>0$
such that $\eta(2t)\leq C\eta(t)$ for all $t>0$.
\end{defi}

A typical element of the class $\mathbb{W}$ is $\eta(t)=t^{1/p}$,
$0<p<\infty$. The functions in the class $\mathbb{W}$ will be used
to define general discrete Lorentz spaces. Occasionally, we will
need to assume a stronger condition on the function $\eta\in
\mathbb{W}$. For $\eta\in\mathbb{W}$ we define the dilation
function
$$M_\eta(s)=\sup_{t>0}\frac{\eta(st)}{\eta(t)}, \; s>0.$$
Since $\eta$ is non-decreasing, $M_\eta(s)\leq 1$ for $0<s\leq 1$.

\begin{defi}\label{d:weightfun_in_W+}
We say that $\eta\in\mathbb{W}_+$ if $\eta\in\mathbb{W}$ and there
exists $s_0\in(0,1)$ for which $M_\eta(s_0)<1$.
\end{defi}

Observe that for $\eta\in\mathbb{W}_+$ and $r>0$,
$\eta^r\in\mathbb{W}_+$. Also, if $\eta\in\mathbb{W}$ and $r>0$,
$t^r\eta(t)\in\mathbb{W}_+$.

\begin{lem}\label{l:sum_etas}
Let $\eta\in \mathbb{W}_+$ and take $s_0$ as in the Definition
\ref{d:weightfun_in_W+}. Then, there exists $C>0$ such that for all
$t>0$
\begin{equation}\label{e:sum_etas}
    \sum_{j=0}^\infty \eta(s_0^jt) \leq C \eta(t).
\end{equation}
\end{lem}

\begin{lem}\label{l:weightfun_C1-diffeomrfsm}
Given $\eta\in\mathbb{W}_+$, there exists $g\in C^1,
g\in\mathbb{W}_+$ such that $g\approx\eta$ and $g'(t)/g(t)\approx
1/t$, $t>0$.
\end{lem}

\vskip0.7cm
\subsection{General Discrete Lorentz Spaces}\label{sS:Gnral_Dscrt_LorentzSpcs}
We will define the discrete Lorentz spaces we will work with.
First, we recall some classical definitions (see e.g. \cite{DeL93}
or \cite{BS}). For a sequence
$\mathbf{s}=\{s_I\}_{I\in\mathcal{D}}\in S$ indexed by the
countable set $\mathcal{D}$, the non-increasing rearrangement of
$\mathbf{s}$ with respect to a measure $\nu$ on $\mathcal{D}$ is
$$\mathbf{s}^\ast_\nu(t)=
    \inf\{\lambda>0: \nu(\{I\in\mathcal{D}:\abs{s_I}>\lambda\})\leq t\}.$$
For $\eta\in\mathbb{W}$, $\nu$ a measure on $\mathcal{D}$, and
$\mu\in(0,\infty]$, the \textbf{discrete Lorentz space}
$\ell^\mu_\eta(\nu)$ is the set of all
$\mathbf{s}=\{s_I\}_{I\in\mathcal{D}}\in S$ such that
\begin{equation}\label{e:def_Gnrl_Dscrt_LorentzSpcs}
\norm{\mathbf{s}}_{\ell^\mu_\eta(\nu)}
    :=\left(\int_0^\infty [\eta(t)\mathbf{s}^\ast_\nu(t)]^\mu \frac{dt}{t}\right)^{1/\mu}<\infty,
\;\;\; 0<\mu<\infty
\end{equation}
and
\begin{equation*}
\norm{\mathbf{s}}_{\ell^\infty_\eta(\nu)}
    := \sup_{t>0} \eta(t)\mathbf{s}^\ast_\nu(t) <\infty.
\end{equation*}

If $\eta(t)=t^{1/p}, 1\leq p<\infty$, then
$\ell^\mu_\eta(\nu)=\ell^{p,\mu}(\nu)$ are the classical
(discrete) Lorentz spaces. For $p=\mu$,
$\ell^{p,p}(\nu)=\ell^p(\nu), 0< p<\infty$, are the spaces of
sequences $\mathbf{s}\in S$ such that
$$\norm{\mathbf{s}}_{\ell^p(\nu)}=\left(\sum_{I\in\mathcal{D}}\abs{s_I}^p\nu(I)\right)^{1/p}.$$
\textbf{Notation}. For $\xi>0$ and $\eta\in\mathbb{W}$,
$\tilde{\eta}(t)=t^\xi\eta(t)\in \mathbb{W}_+$ and
$\ell^\mu_{\tilde{\eta}}(\nu)$ will be denoted by
$\ell^\mu_{\xi,\eta}(\nu)$.

\begin{prop}\label{p:change_notation_LrntzSpcs}
Let $\eta\in\mathbb{W}$ and $\nu$ a measure on $\mathcal{D}$. For
a sequence $\mathbf{s}=\{s_I\}_{I\in\mathcal{D}}\in S$ we have
$$\norm{\mathbf{s}}_{\ell^\infty_\eta(\nu)}\approx
    \sup_{\lambda>0}\lambda\eta(\nu(\{I\in\mathcal{D}:\abs{s_I}>\lambda\})).$$
Moreover, if $0<\mu<\infty$ and $\eta\in\mathbb{W}_+$
$$\norm{\mathbf{s}}_{\ell^\mu_\eta(\nu)}\approx
    \left(\int_0^\infty[\lambda\eta(\nu(\{I\in\mathcal{D}:\abs{s_I}>\lambda\}))]^\mu\frac{d\lambda}{\lambda}\right)^{1/\mu}.$$
\end{prop}

A sequence $\mathbf{u}=\{u_I\}_{I\in\mathcal{D}}\in S$ such that
$u_I > 0$ for all $I\in\mathcal{D}$ will be called a {\bf weight
sequence}.


\begin{defi}\label{d:LrntzSpcs_w_change_of_basis}
Let $\mathbf{u}=\{u_I\}_{I\in\mathcal{D}}$ be a weight sequence and
$\nu$ a positive measure as defined in Subsection
\ref{sS:Rstrct_Nonlin_Apprx_in_Seq_Spcs}. For $0<\mu\leq\infty$ and
$\eta\in\mathbb{W}$ define the space $\ell^\mu_\eta(\mathbf{u},\nu)$
as the set of all sequences $\mathbf{s}=\sum_{I\in\mathcal{D}}
s_I\mathbf{e}_I\in S$ such that
$$ \norm{\mathbf{s}}_{\ell^\mu_\eta(\mathbf{u},\nu)}
    := \norm{\{u_Is_I\}_{I\in\mathcal{D}}}_{\ell^\mu_\eta(\nu)}<\infty.$$
\end{defi}

These spaces will be used in Subsections \ref{sS:J_type_ineqs} and
\ref{sS:B_type_ineqs} to characterize Jackson and Bernstein type
inequalities in the setting of restricted non-linear approximation.
For applications (see Subsections
\ref{sS:Rstrct_Nonlin_Apprx_&_Real_Interpol}-\ref{sS:Appl_Real_Interpol})
we shall take $u_I=\norm{\mathbf{e}_I}_g, I\in\mathcal{D}$, where
$g$ is a quasi-Banach (sequence) lattice.

\begin{lem}\label{l:normalized_ones_bounded}
Let $\mathbf{u}$ and $\nu$ as in Definition
\ref{d:LrntzSpcs_w_change_of_basis} and write
$1_{\Gamma,\mathbf{u}}=\sum_{I\in\Gamma}u_I^{-1}\mathbf{e}_I,
\Gamma\subset\mathcal{D}$ and $\nu(\Gamma)<\infty$.

(a) If $\eta\in\mathbb{W}$,
$\norm{1_{\Gamma,\mathbf{u}}}_{\ell^\infty_\eta(\mathbf{u},\nu)}=\eta(\nu(\Gamma))$.

(b) If $0<\mu<\infty$ and $\eta\in\mathbb{W}$,
$\norm{1_{\Gamma,\mathbf{u}}}_{\ell^\mu_\eta(\mathbf{u},\nu)}\geq
\eta(\nu(\Gamma))$, and if $\eta\in\mathbb{W}_+$,
$\norm{1_{\Gamma,\mathbf{u}}}_{\ell^\mu_\eta(\mathbf{u},\nu)}\approx\eta(\nu(\Gamma))$.
\end{lem}

\vskip0.7cm
\subsection{Jackson type inequalities}\label{sS:J_type_ineqs}
We give equivalent conditions for some Jackson type inequalities
to hold in the setting of restricted non-linear approximation. Our
result generalizes those obtained in \cite{CDH} and \cite{KP06}
for restricted non-linear approximation, as well as those obtained
in \cite{KP04} and \cite{GHN09} for the case $\nu(I)=1$ for all
$I\in\mathcal{D}$ (the counting measure).

\begin{thm}\label{t:Jineq_equiv_for_seq}
Let $(f,\nu)$ be a standard scheme (see Definition
\ref{d:Standard_scheme_Rstrct_NonLin_Apprx})and
 let $\mathbf{u}=\{u_I\}_{I\in\mathcal{D}}$ be a weight sequence. Fix
$\xi>0$ and $\mu\in(0,\infty]$. Then, for any function
$\eta\in\mathbb{W}_+$ the following are equivalent:

1) There exists $C>0$ such that for all $\Gamma\subset\mathcal{D}$
with $\nu(\Gamma)<\infty$
$$\norm{\sum_{I\in\Gamma}\frac{\mathbf{e}_I}{u_I}}_f \leq C \eta(\nu(\Gamma)).$$

2) $\ell^\mu_{\xi,\eta}(\mathbf{u},\nu)\hookrightarrow
\mathcal{A}^\xi_\mu(f,\nu)$.

3)The space $\ell^\mu_{\xi,\eta}(\mathbf{u},\nu)$ satisfies
Jackson's inequality of order $\xi$, that is, there exists $C>0$
such that
$$\sigma_\nu(t,\mathbf{s})_f\leq Ct^{-\xi}\norm{\mathbf{s}}_{\ell^\mu_{\xi,\eta}(\mathbf{u},\nu)}
    , \;\; \text{ for all } \;\;\mathbf{s}\in\ell^\mu_{\xi,\eta}(\mathbf{u},\nu).$$
\end{thm}

Taking $\eta(t)=t^{1/p}, 0<p<\infty$, and
$u_I=\norm{\mathbf{e}_I}_f$ in Theorem \ref{t:Jineq_equiv_for_seq},
condition 1) is called in \cite{KP06} the (upper) p-Temlyakov
property for $f$. In this case,
$\ell^\mu_{\xi,\eta}(\mathbf{u},\nu)=\ell^{q,\mu}(\mathbf{u},\nu)$
with $\frac{1}{q}=\xi+\frac{1}{p}$.

Taking $\nu$ as the counting measure on $\mathcal{D}$ we recover
Theorem 3.6 in \cite{GHN09} from Theorem
\ref{t:Jineq_equiv_for_seq}.

\vskip0.7cm
\subsection{Bernstein type inequalities}\label{sS:B_type_ineqs}
We give equivalent conditions for some Bernstein type inequalities
to hold in the setting of restricted non-linear approximation. This
result generalizes those obtained in \cite{CDH} and \cite{KP06} for
restricted non-linear approximation, as well as those obtained in
\cite{KP04} and \cite{GHN09} for the case $\nu(I)=1$ for all
$I\in\mathcal{D}$ (the counting measure).

We first begin with a representation theorem for the spaces
$\mathcal{A}^\xi_\mu(f,\nu)$. The proof follows that in \cite{P81}
replacing the counting measure by a general positive measure $\nu$.

\begin{prop}\label{p:ReprThm_RstrctApprx}
Let $(f,\nu)$ be a standard scheme (see Definition
\ref{d:Standard_scheme_Rstrct_NonLin_Apprx}). Fix $\xi>0$ and
$\mu\in(0,\infty]$. The following statements are equivalent

i) $\mathbf{s}\in\mathcal{A}^\xi_\mu(f,\nu)$.

ii) There exists $\mathbf{s}_k\in \Sigma_{2^k,\nu}\cap f, k\in\Z$,
such that $\mathbf{s}=\sum_{k=-\infty}^\infty \mathbf{s}_k$ and
$\{2^{k\xi}\norm{\mathbf{s}_k}_f\}_{k\in\Z}\in\ell^\mu(\Z)$.

Moreover,
$$\norm{\mathbf{s}}_{\mathcal{A}^\xi_\mu(f,\nu)}
    \approx \inf\left\{\left[\sum_{k=-\infty}^\infty (2^{k\xi}\norm{\mathbf{s}_k}_f)^\mu\right]^{1/\mu}\right\}$$
where the infimum is taken over all representations of
$\mathbf{s}$ as in ii).
\end{prop}

\begin{thm}\label{t:Bineq_equiv_for_seq}
Let $(f,\nu)$ be a standard scheme (see Definition
\ref{d:Standard_scheme_Rstrct_NonLin_Apprx})and
 let $\mathbf{u}=\{u_I\}_{I\in\mathcal{D}}$ be a weight sequence. Fix $\xi>0$ and
$\mu\in(0,\infty]$. Then, for any function $\eta\in\mathbb{W}$ the
following are equivalent:

1) There exists $C>0$ such that for all $\Gamma\subset\mathcal{D}$
with $\nu(\Gamma)<\infty$,
$$\frac{1}{C}\eta(\nu(\Gamma))\leq \norm{\sum_{I\in\Gamma} \frac{\mathbf{e}_I}{u_I}}_f.$$

2) The space $\ell^\mu_{\xi,\eta}(\mathbf{u},\nu)$ satisfies
Bernstein's inequality of order $\xi$, that is, there exists $C>0$
such that
$$\norm{\mathbf{s}}_{\ell^\mu_{\xi,\eta}(\mathbf{u},\nu)}\leq C t^\xi\norm{\mathbf{s}}_f
    \;\; \text{ for all }\;\; \mathbf{s}\in\Sigma_{t,\nu}\cap f.$$

3) $\mathcal{A}^\xi_\mu(f,\nu)\hookrightarrow
\ell^\mu_{\xi,\eta}(\mathbf{u},\nu)$.
\end{thm}

Taking $\eta(t)=t^{1/p}, 0<p<\infty$, and
$u_I=\norm{\mathbf{e}_I}_f$ in Theorem \ref{t:Bineq_equiv_for_seq},
condition 1) is called in \cite{KP06} the (lower) p-Temlyakov
property for $f$. In this case,
$\ell^\mu_{\xi,\eta}(\mathbf{u},\nu)=\ell^{q,\mu}(\mathbf{u},\nu)$
with $\frac{1}{q}=\xi+\frac{1}{p}$.

Theorem \ref{t:Bineq_equiv_for_seq} generalizes Theorem 5 in
\cite{KP06} for a standard scheme. The proof of this theorem does
not use the theory of real interpolation of quasi-Banach spaces;
we will, however, make use of it to shorten our proof.

Taking $\nu$ as the counting measure on $\mathcal{D}$ we recover
Theorem 4.2 in \cite{GHN09} from Theorem
\ref{t:Bineq_equiv_for_seq}.

\vskip0.7cm
\subsection{Restricted non-linear approximation and real interpolation}\label{sS:Rstrct_Nonlin_Apprx_&_Real_Interpol}
It is well known that $N$-term approximation and real
interpolation are interconnected. If the Jackson and Bernstein's
inequalities hold for $\nu=$ counting measure, $N$-term
approximation spaces are characterized in terms of interpolation
spaces (see e.g. Theorem 3.1 in \cite{DeP88} or Section 9, Chapter
7 in \cite{DeL93}).

As pointed out in \cite{CDH} the above mentioned theory can be
developed in a more general setting. In particular, it can be done
in the frame of the abstract scheme we have introduced in
subsection \ref{sS:Rstrct_Nonlin_Apprx_in_Seq_Spcs}. Below we
state the results we need in this paper. The proofs are
straight-forward modifications of those given in the references
cited in the first paragraph of this section.

\begin{thm}\label{t:J&Bineq_seq->ApprxSpcs_r_InterpolSpcs}
Let $(f,\nu)$ be a standard scheme. Suppose that the quasi-Banach
lattice $g\subset S$ satisfies the Jackson and Bernstein's
inequalities for some $r>0$. Then, for $0<\xi<r$ and
$0<\mu\leq\infty$ we have
$$\mathcal{A}^\xi_\mu(f,\nu)=\left(f,g\right)_{\xi/r,\mu}.$$
\end{thm}

It is not difficult to show that the spaces $\mathcal{A}^r_q(f,\nu),
0<r<\infty, 0<q\leq\infty$, satisfy the Jackson and Bernstein's
inequalities of order $r$, so that by Theorem
\ref{t:J&Bineq_seq->ApprxSpcs_r_InterpolSpcs},
$$\mathcal{A}^\xi_\mu(f,\nu)=\left(f,\mathcal{A}^r_q(f,\nu)\right)_{\xi/r,\mu}$$
for $0<\xi<r$ and $0<\mu\leq\infty$. From here, and using the
reiteration theorem for real interpolation we obtain the following
result that will be used in the proof of Theorem
\ref{t:Bineq_equiv_for_seq}.

\begin{cor}\label{c:Reit_Thm_4_Rstrct_ApprxSpcs}
Let $0<\alpha_0,\alpha_1<\infty, 0<q,q_0,q_1\leq\infty$ and
$0<\theta<1$. Then,
$$\left(\mathcal{A}^{\alpha_0}_{q_0}(f,\nu),
\mathcal{A}^{\alpha_1}_{q_1}(f,\nu)\right)_{\theta,q}
    =\mathcal{A}^\alpha_q(f,\nu), \;\;\; \alpha=(1-\theta)\alpha_0+\theta\alpha_1$$
for a standard scheme $(f,\nu)$.
\end{cor}

\vskip0.7cm
\subsection{Sequence spaces associated with smoothness spaces}
\label{sS:SeqSpcs_Assoctd_w_SmoothSpcs}
A large number of spaces used in Analysis are particular cases of
the Triebel-Lizorkin and Besov spaces.

Given $s\in\R, 0<p<\infty$, and $0<r\leq\infty$, the
Triebel-Lizorkin spaces on $\R^d$ are denoted by
$F^s_{p,r}:=F^s_{p,r}(\R^d)$ where $s$ is a smoothness parameter,
$p$ measures integrability and $r$ measures a refinement of
smoothness. The reader can find the definition of these spaces in
\cite{FJ90, FJW}. Lebesgue spaces $L^p(\R^d)=F^0_{p,2}, 1<p<\infty$,
Hardy spaces $H^p(\R^d)=F^0_{p,2}, 0<p\leq 1$, and Sobolev spaces
$W^s_p(\R^d)=F^s_{p,2}, s>0, 1<p<\infty$, are included in this
collection.

Given $s\in\R, 0<p,r\leq\infty$, the Besov spaces on $\R^d$ are
denoted by $B^s_{p,r}:=B^s_{p,r}(\R^d)$ with an interpretation of
the parameters as in the case of the Triebel-Lizorkin spaces. These
spaces include the Lipschitz classes (see \cite{FJW}).

There are characterizations of $F^s_{p,r}$ and $B^s_{p,r}$ in
terms of sequence spaces. Such characterizations were given first
in \cite{FJ90} using the $\varphi$-transform. Wavelet bases with
appropriate regularity and moment conditions also provide such
characterizations.

A brief description of wavelet bases in $\R^d$ follows. Let
$\mathcal{D}$ be the set of dyadic cubes in $\R^d$ given by
$$Q_{j,k}=2^{-j}([0,1)^d+k), \;j\in\Z, \;k\in\Z^d.$$
A finite collection of functions $\Psi=\{\psi^{(1)},\ldots,
\psi^{(L)}\}\subset L^2(\R^d)$ with $L=2^d-1$ is an (orthonormal)
wavelet family if the set
$$\mathcal{W}:=\{\psi^{(\ell)}_{Q_{j,k}}(x):= 2^{\frac{jd}{2}} \psi^{(\ell)}(2^jx-k)
    :\;Q_{j,k}\in\mathcal{D}, \;\ell=1,2,\ldots,L \}$$
is an orthonormal basis for $L^2(\R^d)$. This is the definition that
appears in \cite{LM86}. The reader can consult properties of
wavelets in \cite{Me90}, \cite{Da92}, \cite{HW96} and \cite{Ma99}.

\begin{defi}\label{d:T-L_assoc_seq_spcs}
Given $s\in\R, 0<p<\infty$ and $0<r\leq\infty$, we let $f^s_{p,r}$
be the space of sequences $\mathbf{s}=\{s_Q\}_{Q\in\mathcal{D}}$
such that
$$\norm{\mathbf{s}}_{f^s_{p,r}}:=\norm{\left[\sum_{Q\in\mathcal{D}} (\abs{Q}^{-s/d+1/r-1/2}\abs{s_Q}\chi^{(r)}_Q(\cdot))^r\right]^{1/r}}_{L^p(\R^d)}<\infty$$
where $\chi^{(r)}_Q(\cdot)= \chi_Q(\cdot)\abs{Q}^{-1/r}$ and
$\chi_Q(\cdot)$ denotes the characteristic function of $Q$.
\end{defi}

\begin{defi}\label{d:B_assoc_seq_spcs}
Given $s\in\R, 0<p,r\leq\infty$, we let $b^s_{p,r}$ be the space
of sequences $\mathbf{s}=\{s_Q\}_{Q\in\mathcal{D}}$ such that
$$\norm{\mathbf{s}}_{b^s_{p,r}}:=\left[\sum_{j\in\Z} \left(
\sum_{\abs{Q}=2^{-jd}}
    (\abs{Q}^{-s/d+1/p-1/2}\abs{s_Q})^p\right)^{r/p}\right]^{1/r}<\infty$$
with the obvious modifications when $p,r=\infty$.
\end{defi}

With appropriate conditions in the elements of a wavelet family
$\Psi=\{\psi^{(1)},\ldots, \psi^{(L)}\}$, $(L=2^d-1)$ it can be
shown that $\mathcal{W}$ is an unconditional basis of $F^s_{p,r}$
or $B^s_{p,r}$ and if $f=\sum_{\ell=1}^L\sum_{Q\in\mathcal{D}}
s^\ell_Q\psi^{(\ell)}_Q$, then
\begin{equation}\label{e:Charactzn_T-L&B_by_seq}
\norm{f}_{F^s_{p,r}}\approx \sum_{\ell=1}^L
\norm{\{s_Q^\ell\}_{Q\in\mathcal{D}}}_{f^s_{p,r}}
    \;\;\; \text{ and }\;\;\;
    \norm{f}_{B^s_{p,r}}\approx \sum_{\ell=1}^L
    \norm{\{s_Q^\ell\}_{Q\in\mathcal{D}}}_{b^s_{p,r}}\,.
\end{equation}
Conditions on $\Psi$ for these equivalences to hold can be found in
\cite{Me90}, \cite{HW96}, \cite{Al}, \cite{Ky}, \cite{LM86}. When a
wavelet family $\Psi$ provides an unconditional basis for
$F^s_{p,r}$ or $B^s_{p,r}$, with equivalences as in
(\ref{e:Charactzn_T-L&B_by_seq}), we shall say that $\Psi$ is
\textbf{admissible} for $F^s_{p,r}$ or $B^s_{p,r}$, respectively.

The equivalences (\ref{e:Charactzn_T-L&B_by_seq}) allow us to work
at the sequence level. We shall drop the sum over $\ell$ since it
only changes the constants in the computations below. The results
proved for sequence spaces $f^s_{p,r}$ or $b^s_{p,r}$ can be
transferred to $F^s_{p,q}$ or $B^s_{p,r}$ by the abstract
transference framework developed in \cite{GH04}.

We notice that the Triebel-Lizorkin and Besov spaces characterized
as in (\ref{e:Charactzn_T-L&B_by_seq}) are called homogeneous,
been often denoted by $\dot{F}^s_{p,r}$ and $\dot{B}^s_{p,r}$. The
non-homogeneous case requires small modifications. Also minor
modifications will allow for the anisotropic spaces as considered
in \cite{GH04}, or the spaces defined by wavelets on bounded
domains. We restrict ourselves to the cases characterized by
(\ref{e:Charactzn_T-L&B_by_seq}).

\vskip0.7cm
\subsection{Restricted approximation for Triebel-Lizorkin sequence spaces}\label{sS:Rstrct_Apprx_T-L_seq-spcs}
As consequence of the theorems developed in Sections
\ref{sS:J_type_ineqs} and \ref{sS:B_type_ineqs} we will obtain
results for restricted approximation in Triebel-Lizorkin sequence
spaces. When coupled with the abstract transference framework
developed in \cite{GH04}, our results generalizes those in
\cite{CDH} and, with minor modifications, those obtained in
\cite{H07}.

\begin{lem}\label{l:bound_sum_cubes}
Let $\Gamma\subset\mathcal{D}$ (not necessarily finite),
$x\in\cup_{Q\in\Gamma}Q$, and $\gamma\not=0$. Define
$$S^\gamma_\Gamma(x)=\sum_{Q\in\Gamma}\abs{Q}^\gamma\chi_Q(x).$$

i) If $\gamma>0$ and there exists $Q^x$, the biggest cube in
$\Gamma$ that contains $x$, then $S^\gamma_\Gamma(x)\approx
\abs{Q^x}^\gamma\chi_{Q^x}(x)=\abs{Q^x}^\gamma$

ii) If $\gamma<0$ and there exists $Q_x$, the smallest cube in
$\Gamma$ that contains $x$, then $S^\gamma_\Gamma(x)\approx
\abs{Q_x}^\gamma\chi_{Q_x}(x)=\abs{Q_x}^\gamma$.
\end{lem}

The smallest cube $Q_x$ from $\Gamma$ that contains
$x\in\cup_{Q\in\Gamma}Q$ has been used by other authors in the
context of non-linear approximation with wavelet basis (see
\cite{HJLY}, \cite{CDH}, \cite{GH04}, \cite{GHM08}). As far as we
know, the biggest cube $Q^x$ from $\Gamma$ that contains
$x\in\cup_{Q\in\Gamma}Q$ has not been used before.

\begin{thm}\label{t:T-L_seq_upper-lower_bounds}
Let $s_1,s_2\in\R$, $0<p_1,p_2<\infty$, $0<q_1,q_2\leq\infty$. For
$\alpha\in \R$ and $\Gamma\subset\mathcal{D}$ define
$\nu_\alpha(\Gamma)=\sum_{Q\in\Gamma}\abs{Q}^\alpha$. Suppose
$\nu_\alpha(\Gamma)<\infty$. Then,
\begin{equation}\label{e:2-10-1}
\norm{\sum_{Q\in\Gamma}\frac{\mathbf{e}_Q}
{\norm{\mathbf{e}_Q}_{f^{s_2}_{p_2,q_2}}}}_{f^{s_1}_{p_1,q_1}}
\approx [\nu_\alpha(\Gamma)]^{1/p_1}
\end{equation}
if and only if $\alpha \neq 1$ and
$\alpha=p_1(\frac{s_2-s_1}{d}-\frac{1}{p_2})+1$  or $\alpha =1,
\frac{s_2-s_1}{d}=\frac{1}{p_2}$ and $p_1=q_1$.
\end{thm}

Theorems \ref{t:Jineq_equiv_for_seq} and \ref{t:Bineq_equiv_for_seq}
with $\eta(t)=t^{1/p_1}$ and
$u_Q=\norm{\mathbf{e}_Q}_{f^{s_2}_{p_2,q_2}}$ together with Theorem
\ref{t:T-L_seq_upper-lower_bounds} show that non-linear
approximation with error measured in $f^{s_1}_{p_1,q_1}$ when the
basis is normalized in $f^{s_2}_{p_2,q_2}$ is related to the use of
the measure $\nu_\alpha(Q)=\abs{Q}^\alpha$, $Q\in\mathcal{D}$,
$\alpha=p_1(\frac{s_2-s_1}{d}-\frac{1}{p_2})$, to control the number
of terms in the approximation. Notice that no role is played by the
second smoothness parameters $q_1,q_2$.

Theorems \ref{t:Jineq_equiv_for_seq} and
\ref{t:Bineq_equiv_for_seq} together with Theorem
\ref{t:T-L_seq_upper-lower_bounds} also allow us to identify the
restricted approximation spaces in the Triebel-Lizorkin setting as
discrete Lorentz spaces.

\begin{cor}\label{c:RstrctApprx(T-L_seq)_r_LrntzSpcs}
Let $s_1,s_2\in\R$, $0<p_1,p_2<\infty$, $0<q_1,q_2\leq\infty$ and
define $\alpha=p_1(\frac{s_2-s_1}{d}-\frac{1}{p_2})+1$. For
$\Gamma\subset\mathcal{D}$ define
$\nu_\alpha(\Gamma)=\sum_{Q\in\Gamma}\abs{Q}^\alpha$. Let $\xi>0$
and $\mu\in(0,\infty]$. If $\alpha\neq 1$,
$$\mathcal{A}^\xi_\mu(f^{s_1}_{p_1,q_1}, \nu_\alpha)
    = \ell^{\tau,\mu}(\mathbf{u},\nu_\alpha),$$
where $\frac{1}{\tau}=\xi+\frac{1}{p_1}$ and
$\mathbf{u}=\{\norm{\mathbf{e}_Q}_{f^{s_2}_{p_2,q_2}}\}_{Q\in\mathcal{D}}$.
If $\alpha =1$ the result holds with $p_1=q_1$.
\end{cor}

For particular values of the parameters, the discrete Lorentz
spaces that appear in Corollary
\ref{c:RstrctApprx(T-L_seq)_r_LrntzSpcs} can be identified as
Besov spaces.

\begin{lem}\label{l:LrntzSpcs_r_BsvSpcs_seq}
Let $s_1,s_2\in\R$, $0<p_1,p_2<\infty$, $0<q_2\leq\infty$ and define
$\alpha=p_1(\frac{s_2-s_1}{d}-\frac{1}{p_2})+1$. For
$\Gamma\subset\mathcal{D}$ define
$\nu_\alpha(\Gamma)=\sum_{Q\in\Gamma}\abs{Q}^\alpha$. Given
$\tau\in(0,\infty)$ we have
$$\ell^{\tau,\tau}(\mathbf{u},\nu_\alpha)=b^\gamma_{\tau,\tau},$$
(with equal quasi-norms) where
$\mathbf{u}=\{\norm{\mathbf{e}_Q}_{f^{s_2}_{p_2,q_2}}\}_{Q\in\mathcal{D}}$
and $\gamma=s_1+d(\frac{1}{\tau}-\frac{1}{p_1})(1-\alpha)$.
\end{lem}

\begin{rem}\label{r:LrntzSpcs_r_BsvSpcs_seq}
If we consider the point of view of \cite{KP06}, then we can only
prove
$\ell^{\tau,\tau}(\mathbf{u},\nu_\alpha)=b^\gamma_{\tau,\tau}\cap
f^{s_1}_{p_1,q_1}$ and the equivalence of quasi-norms holds if we
take $\norm{\cdot}_{b^\gamma_{\tau,\tau}}$ in the right-hand side.

When $\tau<p_1$ and
$\frac{\gamma}{d}-\frac{1}{\tau}=\frac{s_1}{d}-\frac{1}{p_1}$ it is
known that $b^\gamma_{\tau,\tau}\hookrightarrow f^{s_1}_{p_1,q_1}$
(see \cite{FJ88} or \cite{Bu}). This situation occurs when
$\alpha=0$ (the counting measure) but it is not true in our more
general situation.
\end{rem}

The following result identifies certain non-linear approximation
spaces in the restricted setting, when the error is measured in
Triebel-Lizorkin spaces, as Besov spaces. It is obtained as an
easy corollary to Lemma \ref{l:LrntzSpcs_r_BsvSpcs_seq} and
Corollary \ref{c:RstrctApprx(T-L_seq)_r_LrntzSpcs}.

\begin{cor}\label{c:RstrctApprx(T-L_seq)_r_BsvSpcs}
Let $s_1,s_2\in\R, 0<p_1,p_2<\infty, 0<q_1\leq\infty$ and define
$\alpha=p_1(\frac{s_2-s_1}{d}-\frac{1}{p_2})+1$. For
$\Gamma\subset\mathcal{D}$ define
$\nu_\alpha(\Gamma)=\sum_{Q\in\Gamma}\abs{Q}^\alpha$. Given $\xi>0$
define $\tau$ by $\frac{1}{\tau}=\xi+\frac{1}{p_1}$. If $\alpha \neq
1$,
$$\mathcal{A}^\xi_\tau(f^{s_1}_{p_1,q_1}, \nu_\alpha)
    = b^\gamma_{\tau,\tau} \;\;\; \text{(equivalent quasi-norms)},$$
where $\gamma=s_1+d\xi(1-\alpha)$. If $\alpha =1$ the result holds
with $\gamma = s_1$ and $p_1=q_1$.
\end{cor}

\begin{rem}\label{r:RstrctApprx(T-L_seq)_r_BsvSpcs}
If we were to apply Theorem 1 in \cite{KP06} we will obtain
$\mathcal{A}^\xi_\tau(f^{s_1}_{p_1,q_1}, \mathbf{u}, \nu_\alpha)=
b^\gamma_{\tau,\tau}\cap f^{s_1}_{p_1,q_1}$ with equivalence of
quasi-norms, as in Remark \ref{r:LrntzSpcs_r_BsvSpcs_seq}.
\end{rem}

The results obtained in \cite{CDH} for restricted non-linear
approximation with wavelets in the Hardy space $H^p, 0<p<\infty$,
when the wavelets coefficients are restricted to $H^r,
0<r<\infty$, are simple consequences of the above results and the
abstract transference framework developed in \cite{GH04}. To see
this, notice that the sequence spaces associated to $H^p$ and
$H^r$ ($r,p$ as above) are $f^0_{p,2}$ and $f^0_{r,2}$,
respectively.

Thus, for a wavelet basis
$\mathcal{W}=\{\psi^{(\ell)}_Q:Q\in\mathcal{D}, \ell=1,\ldots,
L\}, (L=2^d-1)$ admissible for $H^p$ and $B^\gamma_{\tau,\tau}$,

\begin{equation}\label{e:RA_Hp_r_Bsv_CDH}
\mathcal{A}^\xi_\tau(H^p, \mathcal{W},
\nu_\alpha)=B^\gamma_{\tau,\tau},
\end{equation}
where $\gamma=\frac{dp}{r}\xi$, $\tau$ defined by
$\frac{1}{\tau}=\xi+\frac{1}{p}$, and $\alpha= 1-p/r (\not= 1)$.
This is Corollary 6.3 in \cite {CDH}. Notice that
$\mathcal{A}^\xi_\tau(H^p, \mathcal{W}, \nu_\alpha)$ corresponds to
an approximation space where the wavelet coefficients are normalized
in $H^r$. In the above notation we have emphasize that the
approximation spaces are defined using wavelet basis.

In this situation, The Jackson and Bernstein's inequalities
(Theorems 5.1 and 5.2 in \cite{CDH}) follow from
(\ref{e:RA_Hp_r_Bsv_CDH}) and the fact that the approximation
spaces always satisfy the Jackson and Bernstein's inequalities.

The other situation considered in \cite{CDH} is $B_p:=B^0_{p,p},
0<p<\infty$, when the wavelet coefficients are restricted in $H_r,
0<r<\infty$. In this case, the sequence spaces associated to
$B^0_{p,p}=F^0_{p,p}$ and $H^r$ are $f^0_{p,p}$ and $f^0_{r,2}$,
respectively. Corollary \ref{c:RstrctApprx(T-L_seq)_r_BsvSpcs}
then produces
\begin{equation}\label{e:RA_Bp_r_Bsv_fun}
\mathcal{A}^\xi_\tau(B_p,\mathcal{W},
\nu_\alpha)=B^\gamma_{\tau,\tau},
\end{equation}
where $\gamma=\frac{dp}{r}\xi$, $\tau$ defined by
$\frac{1}{\tau}=\xi+\frac{1}{p}$, with $\nu_\alpha$ and
$\mathcal{W}$ as before. This is more general than Corollary 6.1 in
\cite{CDH} and a comparison with (\ref{e:RA_Hp_r_Bsv_CDH}) proves
immediately a more general version of Theorem 6.3 in \cite{CDH}. Of
course, the Jackson and Bernstein's inequalities of Theorems 5.4 and
5.5 in \cite{CDH} also follow from our results.

To show an example not treated in \cite{CDH} consider the wavelet
orthonormal basis $\mathcal{W}=\{\psi^{(\ell)}_Q:\;Q\in\mathcal{D},
\;\ell=1,\ldots,L \} \;(L=2^d-1)$ admissible for the Sobolev space
$W^s_2, s>0$. We want to measure the error in $W^s_2$ but we
restrict the wavelet coefficients to $L^2$. Since the sequence
spaces associated to $W^s_2$ and $L^2$ are $f^s_{2,2}$ and
$f^0_{2,2}$, Corollary \ref{c:RstrctApprx(T-L_seq)_r_BsvSpcs} and
the abstract framework of \cite{GH04} gives
$$\mathcal{A}^\xi_\tau(W^s_2,\mathcal{W},
\nu_\alpha)=b^\gamma_{\tau,\tau}$$ where $\gamma=s+d\xi(1-\alpha)$,
$\alpha=-2s/d$ and $\tau$ defined by
$\frac{1}{\tau}=\xi+\frac{1}{2}$.

We remark that defining appropriate sequence spaces, a little more
work will show the results proved in \cite{H07} for the anisotropic
case.

\vskip0.7cm
\subsection{Application to Real Interpolation}\label{sS:Appl_Real_Interpol}
Once the restricted approximation spaces for Triebel-Lizorkin
sequence spaces have been identified (see Corollary
\ref{c:RstrctApprx(T-L_seq)_r_BsvSpcs}) we can use Theorem
\ref{t:J&Bineq_seq->ApprxSpcs_r_InterpolSpcs} to obtain results
about real interpolation. This method has been used before (see
\cite{DeP88} or \cite{GH04}). But in the classical case, the
parameters of the spaces interpolated are restricted. With the
theory of restricted approximation we will prove interpolation
results for a much larger set of parameters.

\begin{thm}\label{t:interpol(T-L&Bsv)_r_Bsv_seq}
Let $s\in\R$, $0<p<\infty$, $0<q\leq\infty$. For $0<\tau<p$,
$0<\theta<1$ and $\gamma\not= s$ $(\gamma\in\R)$ we have
$$\left(f^s_{p,q},
b^\gamma_{\tau,\tau}\right)_{\theta,\tau_\theta}=b^{(1-\theta)s+\theta\gamma}_{\tau_\theta,\tau_\theta}
    \;\;\; \text{ with }\;\;\; \frac{1}{\tau_\theta}=\frac{1-\theta}{p}+\frac{\theta}{\tau}.$$
\end{thm}

\begin{rem}\label{r:interpol(T-L&Bsv)_r_Bsv_seq}
Although the Theorem is presented as a result about interpolation
of Triebel-Lizorkin and Besov (sequence) spaces, it is a result
about interpolation of Triebel-Lizorkin sequence spaces, since
$b^\gamma_{\tau,\tau}= f^\gamma_{\tau,\tau}$. Thus, the result can
be stated as
\begin{equation}\label{e:interpol(T-L&T-L)_r_T-L_seq}
\left(f^s_{p,q},
f^\gamma_{\tau,\tau}\right)_{\theta,\tau_\theta}=f^{(1-\theta)s+\theta\gamma}_{\tau_\theta,\tau_\theta}
    \;\;\; \text{ with }\;\;\;
    \frac{1}{\tau_\theta}=\frac{1-\theta}{p}+\frac{\theta}{\tau}.
\end{equation}
\end{rem}

\begin{rem}\label{r:interpol(T-L&Bsv)_r_Bsv_seq2}
Notice that we do not need the restriction
$\frac{\gamma}{d}-\frac{1}{\tau}=\frac{s}{d}-\frac{1}{p}$
characteristic of this type of results when classical non-linear
approximation is used (see e.g. \cite{DeP88} or \cite{GH04}).
\end{rem}

\begin{rem}\label{r:interpol(T-L&Bsv)_r_Bsv_seq3}
By the transference framework designed in \cite{GH04}, the result
of Theorem \ref{t:interpol(T-L&Bsv)_r_Bsv_seq} can be translated
to a result for (homogeneous) Triebel-Lizorkin spaces. The
non-homogeneous case and the case of bounded domains can also be
obtained with minor modifications in the proof.
\end{rem}

\begin{rem}\label{r:interpol(T-L&Bsv)_r_Bsv_seq4}
The merit of Theorem \ref{t:interpol(T-L&Bsv)_r_Bsv_seq} is that
proves results for a large set of parameters using the theory of
approximation. Nevertheless, many (but not all, as far as we know)
have already been proved. One can read from Theorem 3.5 in
\cite{Bu} the following result:
\begin{equation}\label{e:interpol(T-L&T-L)_r_T-L_fun_BUI}
\left(F^{s_0}_{p_0,q_0},
F^{s_1}_{p_1,q_1}\right)_{\theta,p}=F^{s}_{p,p}=B^{s}_{p,p}
\end{equation}
when  $p_i<q_i$, $i=0,1$,
$\frac{1}{p}=\frac{1-\theta}{p_0}+\frac{\theta}{p_1}$,
$s=(1-\theta)s_0+\theta s_1$ and $s_0\not=s_1$. Comparing with
(\ref{e:interpol(T-L&T-L)_r_T-L_seq}) we see that
(\ref{e:interpol(T-L&T-L)_r_T-L_fun_BUI}) has a larger set of
parameters in the second space, while
(\ref{e:interpol(T-L&T-L)_r_T-L_seq}) does not have the restriction
$p< q$ that is required in
(\ref{e:interpol(T-L&T-L)_r_T-L_fun_BUI}). Both of these
shortcomings are due to the methods of the proofs. On the other
hand, Theorem 2.42/1 (page 184) of \cite{Tr78} shows
(\ref{e:interpol(T-L&T-L)_r_T-L_fun_BUI}) without $p_i<q_i$ but
assuming $1<p_i,q_i<\infty$.
\end{rem}

\vskip 1cm
\section{Proofs}\label{S:Proofs}\setcounter{subsection}{3}
\subsection{Weight functions for discrete Lorentz spaces
(proofs)}\label{sS:weightfuns_discrt_LrntzSpcs_Proofs} $\text{ }$

\vskip0.3cm

\textbf{Proof of Lemma \ref{l:sum_etas}}. Let
$\delta:=M_\eta(s_0)<1$. By definition of $M_\eta$ we have
$$1>\delta\geq \frac{\eta(s_0^{j+1}t)}{\eta(s_0^jt)} \;\;\text{
for all }\;\; j=0,1,2,\ldots$$ Therefore,
$$\sum_{j=0}^\infty\eta(s_0^jt)\leq\sum_{j=0}^\infty\delta^j\eta(t)=\eta(t)\frac{1}{1-\delta}.$$
\hfill $\blacksquare$ \vskip .5cm   

\textbf{Proof of Lemma \ref{l:weightfun_C1-diffeomrfsm}}. Define
$g(t)=\int_0^t\frac{\eta(s)}{s}ds$. With $s_0$ as in Definition
\ref{d:weightfun_in_W+}
\begin{eqnarray*}
  g(t)
    &=& \sum_{j=0}^\infty \int_{s_0^{j+1}t}^{s_0^jt}\frac{\eta(s)}{s}ds
        \leq \sum_{j=0}^\infty \eta(s_0^jt)\log(s_0^{-1})
        \leq C\eta(t)\log(s_0^{-1})
\end{eqnarray*}
by Lemma \ref{l:sum_etas} ($C=\frac{1}{1-\delta}$, see the proof
of Lemma \ref{l:sum_etas}). On the other hand
\begin{eqnarray*}
  g(t) &\geq& \int_{t/2}^t\frac{\eta(s)}{s}ds\geq\eta(t/2)\log
  2\geq D\eta(t)\log 2
\end{eqnarray*}
by the doubling property of $\eta$. This shows
\begin{equation}\label{e:Lo-Up_bound_C1-diffeomrfsm}
C_1\eta(t)\leq g(t)\leq C_2\eta(t), \;\; t\in(0,\infty)
\end{equation}
with $0<C_1\leq C_2<\infty$. As an alternative proof one can see
that $\eta$ satisfies the hypotheses of Lemma 1.4 in \cite{KPS92}
(p. 54) to conclude $g\approx \eta$. The function $g$ is clearly
non-decreasing and (\ref{e:Lo-Up_bound_C1-diffeomrfsm}) shows that
$g\in\mathbb{W}$. It is clear that $g\in C^1$ with
$g'(t)=\eta(t)/t$. Thus
$$\frac{g'(t)}{g(t)}=\frac{\eta(t)/t}{\eta(t)}=\frac{1}{t}.$$
It remains to prove that $g\in\mathbb{W}_+$. To prove this,
observe that a function $\eta\in\mathbb{W}$ is an element of
$\mathbb{W}_+$ if and only if
$$i_\eta:=\lim_{t\rightarrow 0^+}\frac{\log M_\eta(t)}{\log t}>0$$
($i_\eta$ is called the lower dilation (Boyd) index of $\eta$ -
see \cite{BS}). Using (\ref{e:Lo-Up_bound_C1-diffeomrfsm})
$$i_g=\lim_{t\rightarrow 0^+} \frac{\log M_g(t)}{\log t}
    \geq \lim_{t\rightarrow 0^+} \frac{\log(\frac{C_1}{C_2}M_\eta(t))}{\log t}
    = \lim_{t\rightarrow 0^+} \frac{\log(M_\eta(t))}{\log t}=i_\eta$$
and, similarly
$$i_g\leq \lim_{t\rightarrow 0^+}\frac{\log(\frac{C_2}{C_1}M_\eta(t))}{\log t}
    = i_\eta.$$
Thus, $i_g=i_\eta>0$ which proves $g\in\mathbb{W}_+$.
\hfill $\blacksquare$ \vskip .5cm   

\vskip0.7cm
\subsection{General discrete Lorentz spaces (proofs)}\label{sS:Gnral_Dscrt_LorentzSpcs_Proofs}$\text{ }$

\vskip0.3cm \textbf{Proof of Proposition
\ref{p:change_notation_LrntzSpcs}.} The case $\mu=\infty$ follows
from part iii) of Proposition 2.2.5 in \cite{CRS}. For
$0<\mu<\infty$, let $w(t)=[\eta(t)]^\mu/t$, $0<t<\infty$. Writing
$\lambda_\nu(t,\mathbf{s})=\nu(\{I\in\mathcal{D}:\abs{s_I}>t\})$ for
the distribution function of $\mathbf{s}$ with respect to the
measure $\nu$ and $\mathrm{W}(s)=\int_0^s w(t)dt$, $0<s<\infty$,
part ii) of Proposition 2.2.5 in \cite{CRS} gives
$$\norm{\mathbf{s}}_{\ell^\mu_\eta(\nu)}
    =\left(\int_0^\infty\mu t^\mu \mathrm{W}(\lambda_\nu(t,\mathbf{s}))\frac{dt}{t}\right)^{1/\mu}.$$
Since $\eta\in \mathbb{W}_+$, $\eta^\mu$ satisfies the hypothesis
of Lemma 1.4 in \cite{KPS92} (p. 54) so that we conclude
$$\mathrm{W}(s)=\int_0^s\frac{\eta(t)^\mu}{t}dt\approx[\eta(s)]^\mu$$
(see also the proof of Lemma \ref{l:weightfun_C1-diffeomrfsm} and
the comment that follows Definition \ref{d:weightfun_in_W+}). This
proves the result.
\hfill $\blacksquare$ \vskip .5cm   

\textbf{Proof of Lemma \ref{l:normalized_ones_bounded}.} (a) Writing
$1_{\Gamma,\mathbf{u}}=\sum_{I\in\mathcal{D}} s_I\mathbf{e}_I$ we
have $s_I=u_I^{-1}$ for all $I\in\Gamma$ and $s_I=0$ if
$I\not\in\Gamma$. Thus, $u_Is_I=1$ for all $I\in\Gamma$ and
$u_Is_I=0$ if $I\not\in\Gamma$. This implies
\begin{equation}\label{e:normalized_ones_bounded_1}
\{u_Is_I\}^\ast_\nu(t)=\left\{
\begin{array}{ll}
    $1$, & \hbox{$0<t<\nu(\Gamma)$} \\
    $0$, & \hbox{$t\geq \nu(\Gamma)$} \\
\end{array}\right\}.
\end{equation}
By Definition \ref{d:LrntzSpcs_w_change_of_basis}
$$\norm{1_{\Gamma,\mathbf{u}}}_{\ell^\infty_\eta(\mathbf{u},\nu)}
    = \sup_{0<t<\nu(\Gamma)}\eta(t)=\eta(\nu(\Gamma)).$$

(b) Using \ref{e:normalized_ones_bounded_1} we have
\begin{eqnarray*}
  \norm{1_{\Gamma,\mathbf{u}}}_{\ell^\mu_\eta(\mathbf{u},\nu)}
    &=& \left(\int_0^{\nu(\Gamma)}[\eta(t)]^\mu\frac{dt}{t}\right)^{1/\mu}
        \geq \left(\int_{\nu(\Gamma)/2}^{\nu(\Gamma)}[\eta(t)]^\mu\frac{dt}{t}\right)^{1/\mu} \\
    &\geq& \eta(\nu(\Gamma)/2)\log 2\geq C\eta(\nu(\Gamma))
\end{eqnarray*}
since $\eta$ is doubling. For the reverse inequality, since
$\eta\in\mathbb{W}_+$, by Proposition
\ref{p:change_notation_LrntzSpcs} we obtain
$$\norm{1_{\Gamma,\mathbf{u}}}_{\ell^\mu_\eta(\mathbf{u},\nu)}\approx
    \left(\int_0^1[\lambda\eta(\nu(\Gamma))]^\mu\frac{d\lambda}{\lambda}\right)^{1/\mu}\approx\eta(\nu(\Gamma)).$$
\hfill $\blacksquare$ \vskip .5cm   

\vskip0.7cm
\subsection{Jackson type inequalities (proofs)}\label{sS:J_type_ineqs_Proofs}
2) $\Rightarrow$ 3) This is immediate since $\cA^\xi_\mu(\nu)
\hookrightarrow \cA^\xi_\infty(\nu)$ and 3) is equivalent to
$\ell^\mu_{\xi,\eta}(\mathbf{u}, \nu) \hookrightarrow
\mathcal{A}^\xi_\infty(f,\nu)$.

3) $\Rightarrow$ 1) Let $0<s_0<1$ be such that $M_\eta(s_0)<1$ as in
the definition of $\eta\in\mathbb{W}_+$. Let
$\Gamma\subset\mathcal{D}$ with $\nu(\Gamma)<\infty$ and write
$1_{\Gamma}:=1_{\Ga,\mathbf{u}}=\sum_{I\in\Ga}\frac{\mathbf{e}_I}{u_I}$.
By Lemma 1 in \cite{KP06} (see also the proof of Theorem 2.1 in
\cite{GH04}), for $\La_0=\Ga\subset\mathcal{D}$ one can find a
subset $\La_1\subset\La_0$ with $\nu(\La_1)\leq s_0\nu(\La_0)$ such
that
$$\norm{1_{\La_0}-1_{\La_1}}_f \approx \si_\nu(s_0\nu(\La_0),1_{\La_0}).$$
We repeat this argument to find nested subsets
$$\Gamma=\Lambda_0\supset\Lambda_1\supset\ldots\supset\Lambda_j\supset\Lambda_{j+1}\supset\ldots$$
such that $\nu(\La_{j+1})\leq s_0\nu(\La_j)$  and
$$\norm{1_{\La_j}-1_{\La_{j+1}}}_f \approx \si_\nu(s_0\nu(\La_j), 1_{\La_j}), \;\;j=0,1,2,\ldots$$
By the $\rho$-power triangle inequality for $f$ we get

\begin{eqnarray*}
  \norm{1_{\Ga}}^\rho_f
    &\leq&  \sum_{j=0}^\infty \norm{1_{\La_j} - 1_{\La_{j+1}}}^\rho_f
        \approx \sum_{j=0}^\infty \si_\nu^\rho(s_0\nu(\La_j), 1_{\La_j}).
\end{eqnarray*}
Using the hypothesis and Lemma \ref{l:normalized_ones_bounded} we
obtain
$$\si_\nu(s_0\nu(\La_j), 1_{\La_j})
    \leq C[s_0\nu(\La_j)]^{-\xi} \norm{1_{\La_j}}_{\ell^\mu_{\xi, \eta}(\nu)}
    \approx \eta(\nu(\La_j))\,. $$
 Concatenating these inequalities we deduce
\begin{eqnarray*}
  \norm{1_{\Ga}}_f
    &\lesssim& \left[ \sum_{j=0}^\infty \eta^\rho (\nu(\La_j)) \right]^{1/\rho}
        \lesssim \left[ \sum_{j=0}^\infty \eta^\rho (s_0^j\nu(\La_0)) \right]^{1/\rho} \\
    &\lesssim& \eta(\nu(\La_0))=\eta(\nu(\Ga))
\end{eqnarray*}
by Lemma \ref{l:sum_etas}, since $\eta$ and $\eta^\rho$ belong to
$\mathbb{W}_+$.

1) $\Rightarrow$ 2) By Lemma \ref{l:weightfun_C1-diffeomrfsm} we may
assume $\eta\in C^1$ and $\eta'(t)/\eta(t)\approx 1/t$, $t>0$. We
start by bounding $\si_\nu(t, \mathbf{s})$ for
$\mathbf{s}\in\ell^\mu_{\xi,\eta}(\mathbf{u},\nu)$. Recall that
$\mathbf{d}:=\{s_Iu_I\}_{I\in\mathcal{D}}\in
\ell^\mu_{\xi,\eta}(\nu)$. Since $\nu(\{I\in\mathcal{D}:
\abs{u_Is_I}> \mathbf{d}^\ast_\nu(t)\})\leq t$ we have

\begin{equation*}
    \si_\nu(t,\mathbf{s}) = \inf_{\mathbf{t}\in \Si_{t,\nu}}
    \norm{\mathbf{s}-\mathbf{t}}_f
        \leq \norm{\sum_{\abs{u_Is_I}\leq \mathbf{d}^\ast_\nu(t)} s_I
        \mathbf{e}_I}_f.
\end{equation*}
For $j=0,1,2,\ldots$ let $\Lambda_j=\{I\in\mathcal{D}:
2^{-j-1}\mathbf{d}^\ast_\nu(t)<\abs{s_Iu_I}\leq
2^{-j}\mathbf{d}^\ast_\nu(t) \}$. The $\rho$-power triangle
inequality and the monotonicity property of $f$, together with the
hypothesis, imply
\begin{eqnarray*}
  [\sigma_\nu(t,\mathbf{s})]^\rho
    &\leq& \sum_{j=0}^\infty \norm{\sum_{I\in\Lambda_j}s_I\mathbf{e}_I}^\rho_f
        = \sum_{j=0}^\infty \norm{\sum_{I\in\Lambda_j}s_Iu_I\frac{\mathbf{e}_I}{u_I}}^\rho_f \\
    &\leq& C \sum_{j=0}^\infty [2^{-j}\mathbf{d}^\ast_\nu(t)]^\rho\eta^\rho(\nu(\Lambda_j)) \\
    &\leq& C \int_0^{\mathbf{d}^\ast_\nu(t)} \lambda^\rho
        [\eta(\nu(\{I\in\mathcal{D}:\abs{s_Iu_I}>\lambda\}))]^\rho\frac{d\lambda}{\lambda}.
\end{eqnarray*}
Applying part 2 of Lemma 2 in \cite{KP06} with
$F(\la)=\frac{\la^\rho}{\rho}$ and $G(\la)=[\eta(\la)]^\rho$ yields
\begin{eqnarray*}
  [\si_\nu(t,\mathbf{s})]^\rho
    &\leq& \frac{1}{\rho} [\mathbf{d}_\nu^\ast(t)\eta(t)]^\rho
        + \frac{1}{\rho}\int_t^{\infty} [\mathbf{d}^\ast_\nu(s)]^\rho d\eta^\rho(s) \\
    &\approx& [\mathbf{d}^\ast_\nu(t)\eta(t)]^\rho
        + \int_t^\infty [\mathbf{d}^\ast_\nu(s)]^\rho \eta^\rho(s)\frac{d s}{s}
\end{eqnarray*}
where we have used $\eta'(s)/\eta(s)\approx 1/s$. Therefore
\begin{equation*}
    \si_\nu(t,\mathbf{s})
        \lesssim  \mathbf{d}^\ast_\nu(t)\eta(t)
        +  \left(\int_t^\infty [\mathbf{d}^\ast_\nu(s)\eta(s)]^\rho \frac{ds}{s}
            \right)^{1/\rho}.
\end{equation*}
Thus,
\begin{eqnarray*}
  \norm{\mathbf{s}}^\mu_{\mathcal{A}^\xi_\mu(\nu)}
    &=& \int_0^\infty[t^\xi\si_\nu(t,\mathbf{s})]^\mu \frac{dt}{t} \\
    &\lesssim&  \int_0^\infty [t^\xi\eta(t)\mathbf{d}_\nu^\ast(t)]^\mu \frac{dt}{t}
        + \int_0^\infty \left[t^\xi\left(\int_t^\infty [\mathbf{d}^\ast_\nu(s)\eta(s)]^\rho \frac{d s}{s}\right)^{1/\rho}\right]^\mu
        \frac{dt}{t}\\
    &:=& I+II.
\end{eqnarray*}
The first term, $I$, is precisely
$\norm{\mathbf{s}}^\mu_{\ell^\mu_{\xi,\eta}(\mathbf{u}, \nu)}$. For
$II$ use Hardy's inequality (see \cite{BS}, p.124) with a $\rho$
such that $\mu/\rho>1$ (notice that this is always possible since if
$f$ satisfies the $\rho$-power triangle inequality it satisfies the
$\rho'$-power triangle inequality for any $0<\rho' \leq \rho$) to
obtain
\begin{eqnarray*}
  II^{1/\mu}
    &=& \left[\int_0^\infty t^{\xi\mu}\left(\int_t^\infty [\mathbf{d}^\ast_\nu(s)\eta(s)]^\rho \frac{d s}{s}\right)^{\mu/\rho}
        \frac{dt}{t}\right]^{1/\mu} \\
    &\leq& \left[ \frac{1}{\rho\xi} \int_0^\infty [\mathbf{d}^\ast_\nu(s) \eta(s)]^\mu s^{\xi \mu}\frac{ds}{s} \right]^{1/\mu} \\
    &=& C \norm{\mathbf{s}}_{\ell^\mu_{\xi,\eta}(\mathbf{u},\nu)},
\end{eqnarray*}
This proves the result.
\hfill $\blacksquare$ \vskip .5cm   

\vskip0.7cm
\subsection{Bernstein type inequalities (proofs)}\label{sS:B_type_ineqs_Proofs}$\text{ }$

\vskip0.3cm

\textbf{Proof of Proposition \ref{p:ReprThm_RstrctApprx}.} i)
$\Rightarrow$ ii) Let $\mathbf{s}\in \mathcal{A}^\xi_\mu(f,\nu)$.
Choose $\mathbf{\varphi}_k\in\Sigma_{2^{k-1},\nu}$ such that
$\norm{\mathbf{s}-\mathbf{\varphi}_k}_f\leq
2\sigma_\nu(2^{k-1},\mathbf{s})$. Let
$\mathbf{s}_{k}=\mathbf{\varphi}_k-\mathbf{\varphi}_{k-1}$, so that
$\mathbf{s}_k\in\Sigma_{2^k,\nu}$. Since
$\mathbf{s}\in\mathcal{A}^\xi_\mu(f,\nu)$ we have
$\sigma_\nu(2^{k-1},\mathbf{s}) \rightarrow 0$ as
$k\rightarrow\infty$; the assumption $f\hookrightarrow S$ implies
$$\lim_{k\rightarrow\infty}\mathbf{\varphi}_k=\mathbf{s}\;\; \text{ in $\mathcal{D}$ (term by term).} $$
On the other hand $\lim_{k\rightarrow-\infty}
\mathbf{\varphi}_k=0$ since $\nu(\text{supp }
\mathbf{\varphi}_k)\rightarrow 0$ as $k\rightarrow-\infty$. Thus,
$$\mathbf{s}=\lim_{k\rightarrow\infty}\mathbf{\varphi}_k=\sum_{k=-\infty}^\infty \mathbf{s}_k.$$
Now,
$$\norm{\mathbf{s}_k}_f^\rho
    \leq \norm{\mathbf{s}-\mathbf{\varphi}_k}_f^\rho + \norm{\mathbf{s}-\mathbf{\varphi}_{k-1}}_f^\rho
    \leq 2\cdot2^\rho[\si_\nu(2^{k-2},\mathbf{s})]^\rho.$$
Therefore,
\begin{equation*}
\sum_{k\in\Z}
    [2^{k\xi}\norm{\mathbf{s}_k}_f]^\mu
    \leq C \sum_{k\in\Z}[2^{k\xi} \si_\nu(2^{k-2},
    \mathbf{s})]^\mu
    \approx \norm{\mathbf{s}}^\mu_{\cA^\xi_\mu(\nu)},
\end{equation*}
by the discrete characterization of the restricted approximation
spaces given in Subsection
\ref{sS:Rstrct_Nonlin_Apprx_in_Seq_Spcs}. It is easy to see that
the result also holds for $\mu=\infty$.

ii) $\Rightarrow$ i) Observe that
$\sum_{k=-\infty}^{\ell-1}\mathbf{s}_k\in\sum_{2^\ell,\nu}$ since
each $\mathbf{s}_k\in\Sigma_{2^k,\nu}$. Take $\rho$ such that
$0<\rho<\mu$ and $\norm{\cdot}_f$ satisfies the $\rho$-power
triangle inequality. We have
\begin{equation*}
[\si_\nu(2^\ell,\mathbf{s})]^\rho
    \leq \norm{\mathbf{s}-\sum_{k=-\infty}^{\ell-1}
    \mathbf{s}_k}^\rho_f
    \leq \sum_{k=\ell}^\infty \norm{\mathbf{s}_k}^\rho_f.
\end{equation*}
With $p=\mu/\rho>1$, (here $0<\mu<\infty$) and $u>0$ such that
$u<\xi\rho$ we have
\begin{eqnarray*}
  \norm{\mathbf{s}}^\mu_{\cA^\xi_\mu(\nu)}
   &\approx& \sum_{\ell=-\infty}^\infty
        [2^{\ell\xi}\si_\nu(2^\ell,\mathbf{s})]^\mu
        \leq \sum_{\ell=-\infty}^\infty [2^{\ell\xi}(\sum_{k=\ell}^\infty \norm{\mathbf{s}_k}^\rho_f)^{1/\rho}]^\mu \\
   &=& \sum_{\ell=-\infty}^\infty 2^{\ell\xi \mu}(\sum_{k=\ell}^\infty 2^{-ku}2^{ku}\norm{\mathbf{s}_k}^\rho_f)^p \\
   &\leq& \sum_{\ell=-\infty}^\infty 2^{\ell\xi \mu} (\sum_{k=\ell}^\infty 2^{-ku p'})^{p/p'} (\sum_{k=l}^\infty 2^{ku p} \norm{\mathbf{s}_k}^\mu_f) \\
   &\approx& \sum_{\ell\in\Z} 2^{\ell(\xi\mu-up)} \left(\sum_{k=\ell}^\infty
   2^{ku p}\norm{\mathbf{s}_k}^\mu_f\right) \\
   &=& \sum_{k\in\Z} 2^{ku p} \norm{\mathbf{s}_k}^\mu_f
   (\sum_{\ell=-\infty}^k 2^{\ell(\xi\mu-up)}) \approx \sum_{k\in\Z} 2^{k\xi\mu} \norm{\mathbf{s}_k}^\mu_f.
\end{eqnarray*}
Since the last expression is finite by hypothesis, we have proved
$\mathbf{s}\in\mathcal{A}^\xi_\mu(f,\nu)$ for $0<\mu<\infty$. For
$\mu=\infty$ we have
\begin{eqnarray*}
  \norm{\mathbf{s}}_{\mathcal{A}^\xi_\infty(\nu)}
    &\approx& \sup_{\ell\in\Z} 2^{\ell\xi}\sigma_\nu(2^\ell,\mathbf{s})
        \leq \sup_{\ell\in\Z} 2^{\ell\xi} (\sum_{k=\ell}^\infty \norm{\mathbf{s}_k}_f^\rho)^{1/\rho}  \\
    &=&  \sup_{\ell\in\Z} 2^{\ell\xi} (\sum_{k=\ell}^\infty 2^{-k\xi\rho}2^{k\xi\rho}\norm{\mathbf{s}_k}_f^\rho)^{1/\rho} \\
    &\leq& (\sup_{k\in\Z} 2^{k\xi}\norm{\mathbf{s}_k}_f) \sup_{\ell\in\Z} 2^{\ell\xi} (\sum_{k=\ell}^\infty 2^{-k\xi\rho})^{1/\rho}
            \approx \sup_{k\in\Z} 2^{k\xi}\norm{\mathbf{s}_k}_f.
\end{eqnarray*}
\hfill $\blacksquare$ \vskip .5cm   

\textbf{Proof of Theorem \ref{t:Bineq_equiv_for_seq}.} 1)
$\Rightarrow$ 2) Let $\mathbf{s}\in\Sigma_{t,\nu}\cap f$ and write
$\mathbf{s}=\{s_I\}_{I\in\mathcal{D}}$. Given
$0<\tau\leq\nu(\Gamma)$ choose
$\Lambda_\tau=\{I\in\Gamma:\abs{u_Is_I}\geq
\mathbf{d}^\ast_\nu(\tau)\}$ where
$\mathbf{d}=\{u_Is_I\}_{I\in\mathcal{D}}$. We have
$\tau\leq\nu(\Lambda_\tau)$ (see (4) in \cite{KP06} or
\cite{BL76}). Applying the hypothesis 1) and the monotonicity of
$f$ we obtain
\begin{eqnarray*}
  \mathbf{d}^\ast_\nu(\tau)\eta(\tau)
    &\leq& \mathbf{d}^\ast_\nu(\tau)\eta(\nu(\Lambda_\tau))
        \lesssim \mathbf{d}^\ast_\nu(\tau)\norm{\sum_{I\in\Lambda_\tau}\frac{\mathbf{e}_I}{u_I}}_f
        \leq \norm{\sum_{I\in\Lambda_\tau} s_I\mathbf{e}_I}_f
            \lesssim \norm{\mathbf{s}}_f.
\end{eqnarray*}
Since $\nu(\{I\in\mathcal{D}:\abs{s_Iu_I}>0\})=\nu(\Gamma)\leq t$
we have $\mathbf{d}^\ast_\nu(\tau)=0$ for $\tau\geq t$. Thus,
\begin{eqnarray*}
  \norm{\mathbf{s}}^\mu_{\ell^\mu_{\xi,\eta}(\mathbf{u},\nu)}
    &=& \int_0^t [\tau^\xi\eta(\tau)\mathbf{d}^\ast_\nu(\tau)]^\mu \frac{d\tau}{\tau}
        \lesssim \norm{\mathbf{s}}^\mu_f \int_0^t \tau^{\xi\mu} \frac{d\tau}{\tau}
        \approx t^{\xi\mu} \norm{\mathbf{s}}^\mu_f.
\end{eqnarray*}
The case $\mu=\infty$ is treated similarly.

2) $\Rightarrow$ 1) Let $1_\Gamma:=
1_{\Gamma,\mathbf{u}}=\sum_{I\in\Gamma}\frac{\mathbf{e}_I}{u_I}$ and
$\nu(\Gamma)=t$, so that $1_\Gamma\in\Sigma_{t,\nu}$. We may assume
$1_\Gamma\in f$, since otherwise the right hand side of 1) is
$\infty$, and the result is trivially true. Hypothesis 2) gives
$$\norm{1_\Gamma}_f\gtrsim t^{-\xi}\norm{1_\Gamma}_{\ell^\mu_{\xi,\eta}(\mathbf{u},\nu)}
    \gtrsim \eta(\nu(\Gamma)),$$
where the last inequality is due to Lemma
\ref{l:normalized_ones_bounded}.

3) $\Rightarrow$ 2) For $\mathbf{s}\in\Sigma_{t,\nu}\cap f$,
$\sigma_\nu(\tau,\mathbf{s})=0$ if $\tau\geq t$. Thus, by 3)
\begin{eqnarray*}
  \norm{\mathbf{s}}_{\ell^\mu_{\xi,\eta}(\mathbf{u},\nu)}
    &\lesssim& \norm{\mathbf{s}}_{\mathcal{A}^\xi_\mu(\nu)}
        = \left(\int_0^t[\tau^\xi\sigma_\nu(\tau,\mathbf{s})]^\mu\frac{d\tau}{\tau}\right)^{1/\mu}\\
    &\lesssim& \norm{\mathbf{s}}_f \left(\int_0^t \tau^{\xi\mu}\frac{d\tau}{\tau}\right)^{1/\mu}
        = t^\xi\norm{\mathbf{s}}_f,
\end{eqnarray*}
where we have used
$\sigma_\nu(\tau,\mathbf{s})\leq\norm{\mathbf{s}}_f$ for all
$\tau>0$.

2) $\Rightarrow$ 3) We have already proved that 1) $\Leftrightarrow$
2); since 1) does not depend on $\xi,\mu$, then 2) holds for all
$\tilde{\xi}>0$ and all $\mu\in(0,\infty]$. For any $\tilde{\xi}>0$
take $\tilde{\rho}$ such that
$\ell^\mu_{\tilde{\xi},\eta}(\mathbf{u},\nu)$ satisfies the
$\tilde{\rho}$-power triangular inequality. By Proposition
\ref{p:ReprThm_RstrctApprx} we can find
$\mathbf{s}_k\in\Sigma_{2^k,\nu}\cap f$, $k\in\Z$, such that
$\mathbf{s}=\sum_{k\in\Z} \mathbf{s}_k$ (in $\mathcal{D}$) and
$$\norm{\mathbf{s}}_{\mathcal{A}^{\tilde{\xi}}_{\tilde{\rho}}(\nu)}
    \approx \left(\sum_{k\in\Z}[2^{k\tilde{\xi}} \norm{\mathbf{s}}_f]^{\tilde{\rho}}\right)^{1/{\tilde{\rho}}}.$$
Applying hypothesis 2) to
$\ell^\mu_{\tilde{\xi},\eta}(\mathbf{u},\nu)$ we obtain
\begin{eqnarray*}
  \norm{\mathbf{s}}^{\tilde{\rho}}_{\ell^\mu_{\tilde{\xi},\eta}(\mathbf{u},\nu)}
    &\lesssim& \sum_{k\in\Z} [\norm{\mathbf{s}_k}_{\ell^\mu_{\tilde{\xi},\eta}(\mathbf{u},\nu)}]^{\tilde{\rho}}
        \lesssim \sum_{k\in\Z} (2^{k\tilde{\xi}}\norm{\mathbf{s}_k}_f)^{\tilde{\rho}}
        \approx \norm{\mathbf{s}}^{\tilde{\rho}}_{\mathcal{A}^{\tilde{\xi}}_{\tilde{\rho}}(\nu)}.
\end{eqnarray*}
This means that for $\mu\in(0,\infty]$ and any $\tilde{\xi}>0$ we
have the continuous inclusion
\begin{equation}\label{e:RASpcs_emmbdd_LrntzSpcs}
\mathcal{A}^{\tilde{\xi}}_{\tilde{\rho}}(\nu) \hookrightarrow
\ell^\mu_{\tilde{\xi},\eta}(\mathbf{u},\nu),
\end{equation}
where $\tilde{\rho}$ is the exponent of the $\tilde{\rho}$-power
triangle inequality for
$\ell^\mu_{\tilde{\xi},\eta}(\mathbf{u},\nu)$. From Corollary
\ref{c:Reit_Thm_4_Rstrct_ApprxSpcs}, for $\xi=(\xi_0+\xi_1)/2$ and
any $\rho\in(0,1]$ we have
$$(\mathcal{A}^{\xi_0}_\rho(\nu), \mathcal{A}^{\xi_1}_\rho(\nu))_{1/2, \mu}=\mathcal{A}_\mu^\xi(\nu).$$
Applying (\ref{e:RASpcs_emmbdd_LrntzSpcs}) with
$\tilde{\xi}=\xi_0$, first, then $\tilde{\xi}=\xi_1$ and
$\rho=\min\{\tilde{\rho}_0,\tilde{\rho}_1\}$ we obtain
$$\mathcal{A}^\xi_\mu(\nu)=(\mathcal{A}^{\xi_0}_\rho(\nu),\mathcal{A}^{\xi_1}_\rho(\nu))_{1/2,\mu}
    \hookrightarrow(\ell^\mu_{\xi_0,\eta}(\mathbf{u},\nu),\ell^\mu_{\xi_1,\eta}(\mathbf{u},\nu))_{1/2,\mu}=\ell^\mu_{\xi,\eta}(\mathbf{u},\nu)$$
where the last equality is a result in real interpolation of
discrete Lorentz spaces that can be found in \cite{M84} (Theorem
3).
\hfill $\blacksquare$ \vskip .5cm   

\vskip0.7cm
\setcounter{subsection}{9}
\subsection{Restricted approximation for Triebel-Lizorkin sequence spaces (proofs)}\label{sS:Rstrct_Apprx_T-L_seq-spcs_(proofs)}$\text{ }$

\vskip0.1cm

\textbf{Proof of Lemma \ref{l:bound_sum_cubes}.} i) It is clear
that $\abs{Q^x}^\gamma\chi_{Q^x}(x)\leq S_\Gamma^\gamma(x)$ since
the right hand side of this inequality contain at least the cube
$Q^x$ (and possibly more). For the reverse inequality we enlarge
the sum defining $S_\Gamma^\gamma(x)$ to include all dyadic cubes
contained in $Q^x$. Therefore,
$$S^\gamma_\Gamma(x)\leq \sum_{Q\subset Q^x:Q\in\mathcal{D}}\abs{Q}^\gamma
    = \sum_{j=0}^\infty (2^{-jd}\abs{Q^x})^\gamma
    = \abs{Q^x}^\gamma\sum_{j=0}^\infty 2^{-jd\gamma}\approx\abs{Q^x}^\gamma$$
since $\gamma>0$.

ii) It is clear that $\abs{Q_x}^\gamma\chi_{Q_x}(x)\leq
S^\gamma_\Gamma(x)$ since the right hand side of this inequality
contains at least the cube $Q_x$ (and possibly more). For the
reverse inequality, we enlarge the sum defining
$S^\gamma_\Gamma(x)$ to include all dyadic cubes containing $Q_x$.
Therefore,
$$S^\gamma_\Gamma(x)\leq \sum_{Q\supset Q_x:Q\in\mathcal{D}}\abs{Q}^\gamma
    = \sum_{j=0}^\infty (2^{jd}\abs{Q_x})^\gamma
    = \abs{Q_x}^\gamma\sum_{j=0}^\infty 2^{jd\gamma}\approx\abs{Q_x}^\gamma$$
since $\gamma<0$.
\hfill $\blacksquare$ \vskip .5cm   

\textbf{Proof of Theorem \ref{t:T-L_seq_upper-lower_bounds}.} We
start by proving (\ref{e:2-10-1}). Write $f_1:=f^{s_1}_{p_1,q_1}$
and $f_2:=f^{s_2}_{p_2,q_2}$ to simplify notation in this proof. By
Definition \ref{d:T-L_assoc_seq_spcs} we have
$\norm{\mathbf{e}_Q}_{f_2}= \abs{Q}^{-s_2/d+1/p_2-1/2}$ and
\begin{eqnarray}\label{e:norm_ones_2diff_T-LSpcs}
\nonumber
\norm{\sum_{Q\in\Gamma}\frac{\mathbf{e}_Q}{\norm{\mathbf{e}_Q}_{f_2}}}_{f_1}
    &=& \left(\int_{\R^d} \left[\sum_{Q\in\Gamma}
        \abs{Q}^{\frac{s_2-s_1}{d}q_1} \abs{Q}^{-q_1/p_2}\chi_Q(x)\right]^{p_1/q_1} dx\right)^{1/p_1} \\
    &=& \left(\int_{\R^d} \left[\sum_{Q\in\Gamma}
        (\abs{Q}^{\frac{\alpha-1}{p_1}} \chi_Q(x))^{q_1}\right]^{p_1/q_1}
        dx\right)^{1/p_1}.
\end{eqnarray}
Consider first the case $\alpha>1$. In this case, since
$\nu_\alpha(\Gamma)<\infty$, the biggest $Q^x$ contained in
$\Gamma$ exists for all $x\in\cup_{Q\in\Gamma}Q$. Applying Lemma
\ref{l:bound_sum_cubes}, part i), first with
$\gamma=\frac{\alpha-1}{p_1}q_1>0$ and then with
$\gamma=\alpha-1>0$, we obtain
$$\left[\sum_{Q\in\Gamma} \abs{Q}^{\frac{\alpha-1}{p_1}q_1}\chi_Q(x)\right]^{p_1/q_1}
    \approx \abs{Q^x}^{\alpha-1}\chi_{Q^x}(x)
    \approx \sum_{Q\in\Gamma} \abs{Q}^{\alpha-1}\chi_Q(x)$$
for all $x\in\cup_{Q\in\Gamma}Q$. From
(\ref{e:norm_ones_2diff_T-LSpcs}) we deduce
\begin{eqnarray*}
  \norm{\sum_{Q\in\Gamma}\frac{\mathbf{e}_Q}{\norm{\mathbf{e}_Q}_{f_2}}}_{f_1}
    &\approx& \left(\int_{\R^d}
        \sum_{Q\in\Gamma} \abs{Q}^{\alpha-1}\chi_Q(x)
        dx\right)^{1/p_1} \\
    &=& \left(\sum_{Q\in\Gamma} \abs{Q}^\alpha\right)^{1/p_1} =
    [\nu_\alpha(\Gamma)]^{1/p_1}.
\end{eqnarray*}

Consider now the case $\alpha<1$. If $\alpha\leq 0$, since
$\nu_\alpha(\Gamma)<\infty$, the smallest cube $Q_x$ contained in
$\Gamma$ exists for all $x\in\cup_{Q\in\Gamma}Q$ (notice that
$\alpha=0$ is  the classical case of counting measure). If
$0<\alpha<1$ we can show that the set $E_\alpha$ of all
$x\in\cup_{Q\in\Gamma}Q$ for which $Q_x$ does not exists has
measure zero. To see this, write
$\mathcal{D}_k=\{Q\in\mathcal{D}:\abs{Q}=2^{-kd}, k\in\Z\}$. Then,
for all $m\geq 0$, $E_\alpha\subset \cup_{k\geq
m}\cup_{Q\in\Gamma\cap\mathcal{D}_k} Q$; therefore
\begin{eqnarray*}
  \abs{E_\alpha}
    &\leq& \sum_{k\geq m}\sum_{Q\in\Gamma\cap\mathcal{D}_k} \abs{Q}
        = \sum_{k\geq m}\sum_{Q\in\Gamma\cap\mathcal{D}_k} \abs{Q}^\alpha\abs{Q}^{1-\alpha} \\
    &\leq& \nu_\alpha(\Gamma)\sum_{k\geq m} 2^{-kd(1-\alpha)}
        \approx \nu_\alpha(\Gamma)2^{-md(1-\alpha)}
\end{eqnarray*}
since $1-\alpha>0$. Letting $m\rightarrow \infty$ we deduce
$\abs{E_\alpha}=0$.

Apply Lemma \ref{l:bound_sum_cubes}, part ii), first with
$\gamma=\frac{\alpha-1}{p_1}q_1<0$ and then with
$\gamma=\alpha-1<0$ to obtain
$$\left[\sum_{Q\in\Gamma} \abs{Q}^{\frac{\alpha-1}{p_1}q_1}\chi_Q(x)\right]^{p_1/q_1}
    \approx \abs{Q_x}^{\alpha-1}\chi_{Q_x}(x)
    \approx \sum_{Q\in\Gamma} \abs{Q}^{\alpha-1}\chi_Q(x)$$
for all $x\in\cup_{Q\in\Gamma}Q$ if $\alpha\leq 0$ and all
$x\in\cup_{Q\in\Gamma}Q\setminus E_\alpha$ if $0<\alpha<1$. In any
case, from (\ref{e:norm_ones_2diff_T-LSpcs}) we deduce
\begin{eqnarray*}
  \norm{\sum_{Q\in\Gamma}\frac{\mathbf{e}_Q}{\norm{\mathbf{e}_Q}_{f_2}}}_{f_1}
    &\approx& \left(\int_{\R^d}
        \sum_{Q\in\Gamma} \abs{Q}^{\alpha-1}\chi_Q(x)
        dx\right)^{1/p_1}
    = \left(\sum_{Q\in\Gamma} \abs{Q}^\alpha\right)^{1/p_1} =
        [\nu_\alpha(\Gamma)]^{1/p_1}.
\end{eqnarray*}

For $\alpha =1$ the set $E_1$ of all $x\in \cup_{Q\in \Gamma} Q$ for
which $Q_x$ exists has also measure zero. Indeed
$$
  \abs{E_1}
    \leq \sum_{k\geq m}\sum_{Q\in\Gamma\cap\mathcal{D}_k} \abs{Q}
        = \sum_{k\geq m}\nu_1 (\Gamma \cap\mathcal{D}_k)
$$
and the last sum tends to zero as $m \to \infty$ since they are the
tails of the convergent sum $\sum_{k\geq 1}\nu_1 (\Gamma
\cap\mathcal{D}_k) \leq \nu_1(\Gamma) <\infty\,.$

Suppose now that (\ref{e:2-10-1}) holds. For $N\in \N$ and $L=2^l$
consider the set $\Gamma_{N,L} = \{[0,L]^d + Lj : j\in \N^d,
0\leq|j|< N\}$ of $N^d$ disjoint dyadic cubes of size length $L$.
For this collection we have
\begin{equation} \label{A}
\nu_\alpha (\Gamma_{N,L}) = \sum_{\Gamma_{N,L}} |Q|^\alpha =
(L^\alpha N)^d\,.
\end{equation}
Also
$$
  \norm{\sum_{Q\in\Gamma_{N,L}}\frac{\mathbf{e}_Q}{\norm{\mathbf{e}_Q}_{f_2}}}_{f_1}
= \left(\int_{\R^d} \left[S_{\Gamma_{N,L}}^\gamma
(x)\right]^{p_1/q_1} dx \right)^{1/p_1}
$$
with $\gamma = (\tfrac{s_2-s_1}{d}-\tfrac{1}{p_2})q_1\,.$ Since
$S_{\Gamma_{N,L}}^\gamma (x) = L^{d\gamma}\sum_{\Gamma_{N,L}}
\chi_Q(x) = L^{d\gamma} \chi_{[0,NL]^d}(x)$ we obtain
\begin{equation} \label{B}
\norm{\sum_{Q\in\Gamma_{N,L}}\frac{\mathbf{e}_Q}{\norm{\mathbf{e}_Q}_{f_2}}}_{f_1}
= L^{d\gamma/q_1}(LN)^{d/p_1} =
L^{d\left(\tfrac{\gamma}{q_1}+\tfrac{1}{p_1}\right)}\,N^{d/p_1}\,.
\end{equation}
Choose $N,N'\in \N$, $L=2^l, L'=2^{l'}$ such that $L^\alpha N =
(L')^\alpha N'$ so that (\ref{A}) implies $\nu_\alpha (\Gamma_{N,L})
= \nu_\alpha (\Gamma_{N',L'})\,.$ By (\ref{e:2-10-1}) and (\ref{B})
we deduce
$$L^{d\left(\tfrac{\gamma}{q_1}+\tfrac{1}{p_1}\right)}\,N^{d/p_1}
\approx
(L')^{d\left(\tfrac{\gamma}{q_1}+\tfrac{1}{p_1}\right)}\,(N')^{d/p_1}\quad
\Leftrightarrow \quad
\left(\frac{L}{L'}\right)^{d\left(\tfrac{\gamma}{q_1}+\tfrac{1-\alpha}{p_1}\right)}\approx
1\,.$$
 This forces $\tfrac{\gamma}{q_1}=\tfrac{\alpha-1}{p_1}$, or
 equivalently $\tfrac{s_2-s_1}{d}- \tfrac{1}{p_2}=\tfrac{\alpha-1}{p_1}$ as desired.

For $\alpha=1$ we still have to prove that $p_1=q_1$. Let $N\in \N$
and $\Gamma_N = \{Q\subset [0,1]^d: 2^{-Nd} < |Q| \leq 1\}\,.$ We
have $\nu_1(\Gamma_N)= N$ and
\begin{equation} \label{C}
\norm{\sum_{Q\in\Gamma_{N}}\frac{\mathbf{e}_Q}{\norm{\mathbf{e}_Q}_{f_2}}}_{f_1}
= \left(\int_{\R^d} \left(\sum_{\Gamma_N} \chi_Q(x)\right)^{p_1/q_1}
dx \right)^{1/p_1} = N^{1/q_1}\,.
\end{equation}
For the same $N\in \N$ take $\tilde \Gamma_N = \{[0,1]^d +
j\overrightarrow{{\bf e}_1}: 0\leq j < N\}$ so that $\nu(\tilde
\Gamma_N) = N$ and
\begin{equation} \label{D}
\norm{\sum_{Q\in\tilde
\Gamma_{N}}\frac{\mathbf{e}_Q}{\norm{\mathbf{e}_Q}_{f_2}}}_{f_1} =
\left(\int_{\R^d}  \chi_{[0,N]\times [0,1]^{d-1}}(x) dx
\right)^{1/p_1} = N^{1/p_1}\,.
\end{equation}
By (\ref{e:2-10-1}) applied to $\Gamma_N$ and $\tilde \Gamma_N$
together with (\ref{C})and (\ref{D}) we obtain $N^{1/q_1} \approx
N^{1/p_1}$. This forces $p_1=q_1$ as we wanted.
\hfill $\blacksquare$ \vskip .5cm   

\textbf{Proof of Corollary
\ref{c:RstrctApprx(T-L_seq)_r_LrntzSpcs}.} Apply Theorems
\ref{t:Jineq_equiv_for_seq} and \ref{t:Bineq_equiv_for_seq} to
$f=f^{s_1}_{p_1,q_1}$, ${\bf u}$ and $\nu_\alpha$, as given in the
statement of the corollary, and with $\eta(t)=t^{1/p_1}$.
\hfill $\blacksquare$ \vskip .5cm   

\textbf{Proof of Lemma \ref{l:LrntzSpcs_r_BsvSpcs_seq}.} Let
$f_2:=f^{s_2}_{p_2,q_2}$ to simplify notation.  Since
$\norm{\mathbf{e}_Q}_{f_2}=\abs{Q}^{-s_2/d+1/p_2-1/2}=\abs{Q}^{-\gamma/d+(1-\alpha)/\tau-1/2}$
 for
$\mathbf{s}=\sum_{Q\in\mathcal{D}}s_Q\mathbf{e}_Q$ we have
\begin{eqnarray*}
  \norm{\mathbf{s}}^\tau_{\ell^{\tau,\tau}(\mathbf{u},\nu_\alpha)}
    &=& \norm{\{\norm{s_Q\mathbf{e}_Q}_{f_2}\}_{Q\in\mathcal{D}}}^\tau_{\ell^{\tau,\tau}
    (\nu_\alpha)}
        = \sum_{Q\in\mathcal{D}} \norm{s_Q\mathbf{e}_Q}^\tau_{f_2}\abs{Q}^\alpha \\
    &=& \sum_{Q\in\mathcal{D}} (\abs{s_Q}\abs{Q}^{-\gamma/d+(1-\alpha)/\tau-1/2}
    \abs{Q}^{\alpha/\tau})^\tau \\
    &=& \sum_{Q\in\mathcal{D}}
        (\abs{s_Q}\abs{Q}^{-\gamma/d+1/\tau-1/2})^\tau =
        \norm{\mathbf{s}}^\tau_{b^\gamma_{\tau,\tau}}.
\end{eqnarray*}
\hfill $\blacksquare$ \vskip .5cm   

\vskip0.7cm
\subsection{Application to real interpolation (proofs)}\label{sS:Appl_Real_Interpol_(proofs)}$\text{ }$

\vskip0.3cm

\textbf{Proof of Theorem \ref{t:interpol(T-L&Bsv)_r_Bsv_seq}.} Write
$\xi=1/\tau-1/p>0$ and choose $\alpha\not=1$ such that
$\gamma=s+(1-\alpha)\xi d$ (\emph{i.e.}
$\alpha=1-\frac{\gamma-s}{\xi d}$), which is possible since
$\gamma\not= s$. Once $\alpha$ is chosen, take $s_2\in\R$ in such a
way that $\alpha=p(\frac{s_2-s}{d})$ (\emph{i.e.}
$s_2=s+\tfrac{\alpha d}{p}$). Theorem
\ref{t:T-L_seq_upper-lower_bounds} shows that $f^s_{p,q}$ satisfies
1) of Theorems \ref{t:Jineq_equiv_for_seq} and
\ref{t:Bineq_equiv_for_seq} with $\eta(t)=t^{1/p}$, for the
"normalization" space $f_2:=f^{s_2}_{p,p}$ and $\nu_\alpha$ (notice
that $\alpha=p(\frac{s_2-s}{d})$ is the condition required in Lemma
\ref{l:LrntzSpcs_r_BsvSpcs_seq}). Thus, the space
$\ell^{\tau,\mu}(\mathbf{u},\nu_\alpha)$ satisfies the Jackson and
Bernstein's inequalities of order $\xi=1/\tau-1/p>0$, where
$\mathbf{u}=\{\norm{\mathbf{e}_Q}_{f_2}\}_{Q\in\mathcal{D}}$. Taking
$\mu=\tau$, Lemma \ref{l:LrntzSpcs_r_BsvSpcs_seq} shows that
$b^\gamma_{\tau,\tau}$ satisfies the Jackson and Bernstein's
inequalities of order $\xi=1/\tau-1/p>0$, since
$\gamma=s+(1-\alpha)\xi d$ (the required condition).

By Theorem \ref{t:J&Bineq_seq->ApprxSpcs_r_InterpolSpcs}, for
$0<\theta<1$,
$$(f^s_{p,q}, b^\gamma_{\tau,\tau})_{\theta,\tau_\theta}
    =\mathcal{A}^{\theta\xi}_{\tau_\theta}(f^s_{p,q}, \nu_\alpha).$$
Since $\frac{1}{\tau_\theta}=
\frac{(1-\theta)}{p}+\frac{\theta}{\tau}=\theta(\frac{1}{\tau}-\frac{1}{p})+\frac{1}{p}=\theta\xi+\frac{1}{p}$,
Corollary \ref{c:RstrctApprx(T-L_seq)_r_BsvSpcs} gives
$$\mathcal{A}^{\theta\xi}_{\tau_\theta}(f^s_{p,q}, \nu_\alpha)
    =b^{\tilde{\gamma}}_{\tau_\theta,\tau_\theta}$$
with $\tilde{\gamma}=
s+d\theta\xi(1-\alpha)=s+d\theta\xi\frac{(\gamma-s)}{\xi d}
=s+\theta(\gamma-s)=(1-\theta)s+\theta\gamma$, which proves the
result.
\hfill $\blacksquare$ \vskip .5cm   

{\bf Acknowledgements}. We thank Gustavo Garrig\'os for reading a
first manuscript of this work and for his suggestions to improve the
presentation.


\end{document}